\documentclass[12pt]{amsart}

\usepackage{amssymb}
\usepackage{graphicx}
\usepackage{hyperref}
\usepackage{ifthen}
\usepackage{mathtools}
\usepackage{pict2e}
\usepackage{xargs}
\usepackage{xspace}

\newcommand{\definedterm}[1]{\emph{#1}}

\newcommand{\action}{\curvearrowright}
\newcommand{\alignedmaps}[1]{\mathrm{Al}_{#1}}
\newcommand{\bbF}{\mathbb{F}}
\newcommand{\bbZ}{\mathbb{Z}}
\newcommand{\brokenpair}[2]{(#1,$ $#2)}
\newcommand{\calA}{\mathcal{A}}
\newcommand{\calB}{\mathcal{B}}
\newcommand{\calF}{\mathcal{F}}
\newcommand{\calI}{\mathcal{I}}
\newcommand{\calJ}{\mathcal{J}}
\newcommand{\calN}{\mathcal{N}}
\newcommand{\calP}{\mathcal{P}}
\newcommand{\calR}{\mathcal{R}}
\newcommand{\calU}{\mathcal{U}}
\newcommand{\calV}{\mathcal{V}}
\newcommand{\calW}{\mathcal{W}}

\newcommandx{\CantorCantorspace}[2][1 =, 2 =]{
  \ifthenelse{\equal{#2}{}}
    {\ifthenelse{\equal{#1}{}}
        {\functions{\N}{(\Cantorspace)}}
        {\functions{#1}{(\Cantorspace)}}}
    {\functions{#2}{(\functions{#1}{2})}}
}
\newcommand{\Cantorspace}[1][]{
  \ifthenelse{\equal{#1}{}}{\functions{\N}{2}}{\functions{#1}{2}}
}
\newcommand{\cardinality}[1]{|#1|}
\newcommand{\case}[2]{\emph{Case #1: #2}}
\newcommand{\compatible}[2]{\mathrm{Comp}(#1, #2)}
\newcommand{\completeer}[1]{I(#1)}
\newcommand{\composition}{\circ}
\newcommandx{\concatenation}[2][1 = undefined, 2 = undefined]{
  \ifthenelse{\equal{#1}{undefined}}{{}\smallfrown}{
    \ifthenelse{\equal{#2}{undefined}}{\bigoplus #1}{\bigoplus_{#1} #2}
  }
}
\newcommand{\configurations}[1]{\mathrm{Conf}_{#1}}
\newcommand{\continuousfunctions}[2]{C(#1, #2)}
\newcommand{\defeq}{\coloneqq}
\newcommand{\diagonal}[1]{\Delta(#1)}
\newcommand{\diam}{\mathrm{diam}}
\newcommandx{\disjointunion}[2][1 =, 2 =]{
  \ifthenelse{\equal{#1}{}}{\amalg}{
    \ifthenelse{\equal{#2}{}}{\coprod #1}{\coprod_{#1} #2}
  }
}

\newcommand{\doublehorizontalsection}[2]{#1^{(#2)}}
\newcommand{\Eone}{\mathbb{E}_1}

\newcommand{\equivalenceclass}[2]{[#1]_{#2}}
\newcommand{\evaluation}{\mathrm{eval}}
\newcommand{\existnonmeagerlymany}{\exists^*}
\newcommand{\extendedby}{\sqsubseteq}
\newcommand{\extensions}[2][]{
  \ifthenelse{\equal{#1}{}}{\calN_{#2}}{\calN_{#2} \intersection #1}
}
\newcommand{\forcomeagerlymany}{\forall^*}
\newcommand{\from}{\colon}
\newcommand{\Fsigma}{$F_\sigma$\xspace}
\newcommandx{\functions}[3][3 =]{
  \ifthenelse{\equal{#3}{}}{#2^{#1}}{#2^{#1}_{#3}}
}
\newcommandx{\Fzero}[3][2 =, 3 =]{
  \ifthenelse{\equal{#2}{}}
    {\bbF_{#1}}
    {\bbF_{#1}(\CantorCantorspace[#2][#3])}
}
\newcommand{\generateder}[1]{\langle #1 \rangle}

\newcommand{\horizontalconcatenation}{\concatenation}
\newcommand{\horizontalsection}[2]{#1^{#2}}
\newcommand{\id}[1]{\mathrm{id}_{#1}}

\newcommandx{\identityfunction}[2][2 =]{
  \ifthenelse{\equal{#2}{}}
    {\mathrm{id} \from #1 \to #1}
    {\mathrm{id} \from #1 \to #2}
}
\newcommand{\image}[2]{#1(#2)}
\newcommand{\intersectingrelation}[1]{#1^{\cap}}
\newcommandx{\intersection}[2][1 =, 2 =]{
  \ifthenelse{\equal{#1}{}}{\cap}{
    \ifthenelse{\equal{#2}{}}{\bigcap #1}{{\bigcap_{#1} #2}}
  }
}
\newcommand{\jump}[1]{#1^+}
\newcommand{\mathand}{\text{ and }}
\newcommand{\N}{\mathbb{N}}
\newcommand{\orbitequivalencerelation}[2]{E_{#1}^{#2}}
\newcommand{\pair}[2]{(#1, #2)}

\renewcommand{\phi}{\varphi}
\newcommandx{\Piclass}[2][1=,2=]{
  \ifthenelse{\equal{#2}{}}{\mathbf{\Pi}_{#1}}{\mathbf{\Pi}^{#1}_{#2}}
}
\newcommand{\preimage}[2]{#1^{-1}(#2)}
\newcommandx{\product}[2][1 =, 2 =]{
  \ifthenelse{\equal{#1}{}}{\times}{
    \ifthenelse{\equal{#2}{}}{\prod #1}{{\prod_{#1} #2}}
  }
}
\newcommandx{\projection}[2][1 =, 2 =]{
  \ifthenelse{\equal{#1}{}}{\mathrm{proj}}{
    \ifthenelse{\equal{#2}{}}{\projection_{#1}}{
      \image{\projection[#1]}{#2}
    }
  }
}
\newcommand{\restrictedalignedmaps}[2]{\calU_{#1}^{#2}}
\newcommandx{\restrictedcontinuousfunctions}[4][3 = undefined,
  4 = undefined]{
    \ifthenelse{\equal{#3}{undefined}}
      {M(#1, #2)}
      {M_{\continuousfunctions{#3}{#4}}(#1, #2)}
  }
\newcommand{\saturation}[2]{[#1]_{#2}}
\newcommand{\sectionpower}[2]{#1^{\rightarrow #2}}
\newcommandx{\sequence}[2][2 = undefined]{
  \ifthenelse{\equal{#2}{undefined}}{(#1)}{
    (#1)_{#2}
  }
}
\newcommandx{\set}[2][2 = undefined]{
  \ifthenelse{\equal{#2}{undefined}}{\{ #1 \}}{
    \{ #1 \suchthat #2 \}
  }
}
\newcommand{\setcomplement}[1]{\twiddle #1}
\newcommand{\setinterval}[3][]{[#2, #3]_{#1}}
\newcommandx{\sets}[4][3 =, 4 =]{
  \ifthenelse{\equal{#4}{}}
    {\ifthenelse{\equal{#3}{}}{[#2]^{#1}}{[#2]^{#1}_{#3}}}
    {[#2]^{#1}_{#3, #4}}
}
\newcommandx{\Sigmaclass}[2][1=,2=]{
  \ifthenelse{\equal{#2}{}}{\mathbf{\Sigma}_{#1}}{\mathbf{\Sigma}^{#1}_{#2}}
}
\newcommand{\suchthat}{\mid}
\newcommand{\triple}[3]{(#1, #2, #3)}
\newcommand{\twiddle}{\raisebox{1pt}{\scalebox{.75}{$\mathord{\sim}$}}}
\newcommand{\uniformmetric}[1]{#1_\infty}
\newcommandx{\union}[2][1 =, 2 =]{
  \ifthenelse{\equal{#1}{}}{\cup}{
    \ifthenelse{\equal{#2}{}}{\bigcup #1}{{\bigcup_{#1} #2}}
  }
}
\newcommand{\verticalconcatenation}{\oplus}
\newcommand{\verticalsection}[2]{#1_{#2}}
\newcommand{\Z}{\bbZ}

\newcommand{\Baire}{Baire\xspace}
\newcommand{\Borel}{Bor\-el\xspace}
\newcommand{\Caicedo}{Cai\-ce\-do\xspace}
\newcommand{\Cantor}{Can\-tor\xspace}
\newcommand{\Cauchy}{Cau\-chy\xspace}
\newcommand{\Clemens}{Clem\-ens\xspace}
\newcommand{\Conley}{Con\-ley\xspace}
\newcommand{\Feldman}{Feld\-man\xspace}
\newcommand{\Friedman}{Fried\-man\xspace}
\newcommand{\Jankov}{Jan\-kov\xspace}
\newcommand{\Kechris}{Kech\-ris\xspace}
\newcommand{\Kuratowski}{Kur\-at\-ow\-ski\xspace}
\newcommand{\Louveau}{Lou\-veau\xspace}
\newcommand{\Lusin}{Lu\-sin\xspace}
\newcommand{\Miller}{Mil\-ler\xspace}
\newcommand{\Moore}{Moore\xspace}
\newcommand{\Mycielski}{My\-ciel\-ski\xspace}
\newcommand{\Novikov}{No\-vik\-ov\xspace}
\newcommand{\Polish}{Po\-lish\xspace}
\newcommand{\Silver}{Sil\-ver\xspace}

\newcommand{\Stanley}{Stan\-ley\xspace}
\newcommand{\Ulam}{U\-lam\xspace}
\newcommand{\vonNeumann}{von Neu\-mann\xspace}

\newenvironment{claimproof}{
  
  \begin{proof}
}{\end{proof}}

\newenvironment{corollaryproof}{
  
  \begin{proof}
}{\end{proof}}

\newenvironment{lemmaproof}{
  
  \begin{proof}
}{\end{proof}}

\newenvironment{propositionproof}{
  
  \begin{proof}
}{\end{proof}}

\newenvironment{theoremproof}{
  
  \begin{proof}
}{\end{proof}}

\newtheorem{claim}{Claim}[section]
\newtheorem{corollary}[claim]{Corollary}
\newtheorem{lemma}[claim]{Lemma}
\newtheorem{proposition}[claim]{Proposition}
\newtheorem{theorem}[claim]{Theorem}

\newtheorem{introtheorem}{Theorem}

\theoremstyle{definition}
\newtheorem*{acknowledgements}{Acknowledgements}
\newtheorem{remark}[claim]{Remark}

\begin{document}

\begin{abstract}
  We show that if an equivalence relation $E$ on a \Polish space is a
  countable union of smooth \Borel subequivalence relations, then there is
  either a \Borel reduction of $E$ to a countable \Borel equivalence relation
  on a \Polish space or a continuous embedding of $\Eone$ into $E$. We
  also establish related results concerning countable unions of more general
  \Borel equivalence relations.
\end{abstract}

\author{N. de Rancourt}

\address{
  N. de Rancourt \\
  Faculty of Mathematics and Physics \\
  Department of Mathematical Analysis \\
  Sokolovsk\'a 49/83 \\
  186 75 Praha 8 \\
  Czech Republic
 }

\email{rancourt@karlin.mff.cuni.cz}

\urladdr{
  \url{https://www2.karlin.mff.cuni.cz/~rancourt/}
}

\author{B.D. Miller}

\address{
  B.D. Miller \\
  Faculty of Mathematics \\
  University of Vienna \\
  Oskar Morgenstern Platz 1 \\
  1090 Wien \\
  Austria
 }

\email{benjamin.miller@univie.ac.at}

\urladdr{
  \url{https://homepage.univie.ac.at/benjamin.miller/}
}

\thanks{The authors were supported, in part, by FWF Grant P29999.}
  


\title[A dichotomy for countable unions]{A dichotomy for countable unions
  of smooth Borel equivalence relations}

\maketitle

\section*{Introduction}

A \definedterm{homomorphism} from a binary relation $R$ on a set $X$ to
a binary relation $S$ on a set $Y$ is a function $\phi \from X \to Y$ with the
property that $\image{(\phi \times \phi)}{R} \subseteq S$. More generally, a
\definedterm{homomorphism} from a sequence $\sequence{R_i}[i \in I]$ of
binary relations on $X$ to a sequence $\sequence{S_i}[i \in I]$ of binary
relations on $Y$ is a function $\phi \from X \to Y$ that is a homomorphisms
from $R_i$ to $S_i$ for all $i \in I$. A \definedterm{cohomomorphism} is a
homomorphism of the corresponding complements, a \definedterm
{reduction} is a homomorphism that is also a cohomorphism, and an
\definedterm{embedding} is an injective reduction.

A \definedterm{\Polish space} is a second countable topological space that
admits a compatible complete metric. A subset of a topological space is
\definedterm{\Borel} if it is in the smallest $\sigma$-algebra containing the
open sets. A function between topological spaces is \definedterm{\Borel} if
preimages of open sets are \Borel.

Following the usual abuse of language, we say that an equivalence relation
is \definedterm{countable} if each of its classes is countable. A \Borel
equivalence relation $E$ on a \Polish space is \definedterm{smooth} if
there is a \Borel reduction of $E$ to equality on a \Polish space,
\definedterm{hypersmooth} if it is the union of an increasing sequence
$\sequence{E_n}[n \in \N]$ of smooth \Borel subequivalence relations, and
\definedterm{$\sigma$-smooth} if it is a countable union of smooth \Borel
subequivalence relations. A well-known example of a hypersmooth \Borel
equivalence relation on $\functions{\N}{(\Cantorspace)}$ that is not \Borel
reducible to a countable \Borel equivalence relation on a \Polish space is
given by $c \mathrel{\Eone} d \iff \exists n \in \N \forall m \ge n \ c(m) =
d(m)$.

Our primary result here is the following analog of the \Kechris--\Louveau
dichotomy for hypersmooth \Borel equivalence relations on \Polish spaces
(see \cite[Theorem 1]{KechrisLouveau} or Theorem \ref{KechrisLouveau}):

\begin{introtheorem} \label{introduction:main}
  Suppose that $E$ is a $\sigma$-smooth \Borel equivalence relation on a
  \Polish space. Then exactly one of the following holds:
  \begin{enumerate}
    \item There is a \Borel reduction of $E$ to a countable \Borel equivalence
      relation on a \Polish space.
    \item There is a continuous embedding of $\Eone$ into $E$.
  \end{enumerate}
\end{introtheorem}

A \definedterm{$\sigma$-ideal} on a set $X$ is a family $\calI$ of subsets
of $X$ that is closed under containment and countable unions. When $X$
is a \Polish space, we say that such a $\sigma$-ideal is \definedterm
{weakly ccc-on-\Borel} if there is no uncountable family of pairwise disjoint
\Borel subsets of $X$ that are not in $\calI$. Given sets $X$ and $Y$, the
\definedterm{horizontal section} of a set $R \subseteq X \times Y$ at a
point $y$ of $Y$ is given by $\horizontalsection{R}{y} \defeq \set{x \in X}[x
\mathrel{R} y]$, and the \definedterm{vertical section} of $R$ at point $x$
of $X$ is given by $\verticalsection{R}{x} \defeq \set{y \in Y}[x \mathrel{R}
y]$. An assignment $x \mapsto \calI_x$, sending each point of $X$ to a
$\sigma$-ideal on $X$, is \definedterm{\Borel-on-\Borel} if $\set{x \in X}
[\verticalsection{R}{x} \in \calI_x]$ is \Borel for all \Borel sets $R \subseteq
X \times X$, and \definedterm{strongly \Borel-on-\Borel} if $\set{\pair{x}{y}
\in X \times Y}[\verticalsection{R}{\pair{x}{y}} \in \calI_x]$ is \Borel for all
\Polish spaces $Y$ and \Borel sets $R \subseteq (X \times Y) \times X$. A
\Borel equivalence relation $E$ on a \Polish space $X$ is \definedterm
{idealistic} if there is an $E$-invariant \Borel-on-\Borel assignment $x
\mapsto \calI_x$ sending each point in $X$ to a $\sigma$-ideal on $X$ for
which $\equivalenceclass{x}{E} \notin \calI_x$. We say that $E$ is
\definedterm{ccc idealistic} if each $\calI_x$ can be taken to be weakly
ccc-on-\Borel, \definedterm{strongly idealistic} if the assignment $x \mapsto
\calI_x$ can be taken to be strongly \Borel-on-\Borel, and \definedterm
{strongly ccc idealistic} if $x \mapsto \calI_x$ can be taken to be a strongly
\Borel-on-\Borel assignment of weakly ccc-on-\Borel $\sigma$-ideals.

Recall that the \definedterm{orbit equivalence relation} induced by a group
action $\Gamma \action X$ is the equivalence relation
$\orbitequivalencerelation{\Gamma}{X}$ on $X$ given by $x \mathrel
{\orbitequivalencerelation{\Gamma}{X}} y \iff \exists \gamma \in \Gamma
\ x = \gamma \cdot y$. The \Feldman--\Moore theorem ensures that every
countable \Borel equivalence relation on a \Polish space is the orbit
equivalence relation induced by a \Borel action of a countable discrete
group (see \cite[Theorem 1]{FeldmanMoore}), and the proof of \cite[\S1.II.i]
{Kechris:LocallyCompact} shows that every \Borel orbit equivalence
relation induced by a \Borel action of a \Polish group on a \Polish space is
strongly ccc idealistic. By \cite[Theorem 4.1]{KechrisLouveau}, the
equivalence relation $\Eone$ is not \Borel reducible to a ccc-idealistic
\Borel equivalence relation on a \Polish space.

Much as the \Kechris--\Louveau dichotomy can be used to show that a
hypersmooth \Borel equivalence relation on a \Polish space is \Borel
reducible to a ccc-idealistic \Borel equivalence relation on a \Polish space
if and only if it is \Borel reducible to the orbit equivalence relation induced
by a \Borel action of $\Z$ on a \Polish space, Theorem \ref
{introduction:main} yields:

\begin{introtheorem} \label{introduction:action}
  Suppose that $E$ is a $\sigma$-smooth \Borel equivalence relation on a
  \Polish space. Then the following are equivalent:
  \begin{enumerate}
    \item There is a \Borel reduction of $E$ to a ccc-idealistic \Borel
      equivalence relation on a \Polish space.
    \item There is a \Borel reduction of $E$ to the orbit equivalence
      relation induced by a \Borel action of a countable discrete group on
      a \Polish space.
  \end{enumerate} 
\end{introtheorem}

A subset of a topological space is \definedterm{\Fsigma} if it is a countable
union of closed sets, and a binary relation $R$ on a \Polish space $X$ is
\definedterm{potentially \Fsigma} if there is a \Polish topology on $X$,
generating the same \Borel sets as the given topology, with respect to
which $R$ is \Fsigma. Standard change of topology results and the
\Lusin--\Novikov uniformization theorem (see, for example, \cite[\S13 and
Theorem 18.10]{Kechris}) easily imply that countable \Borel equivalence
relations on \Polish spaces are potentially \Fsigma.

\Kechris--\Louveau have asked whether a \Borel equivalence relation $E$
on a \Polish space is \Borel reducible to a ccc-idealistic \Borel equivalence
relation on a \Polish space if and only if there is no continuous embedding
of $\Eone$ into $E$. Much as the \Kechris--\Louveau dichotomy yields a
positive answer to this question in the hypersmooth case, Theorem \ref
{introduction:main} yields the extension to the $\sigma$-smooth case, and
the underlying argument can be used to obtain a further generalization:

\begin{introtheorem} \label{introduction:idealistic}
  Suppose that $E$ is an equivalence relation on a \Polish space that is a
  countable union of subequivalence relations that are \Borel reducible to
  strongly-ccc-idealistic potentially-\Fsigma equivalence relations on \Polish
  spaces. Then exactly one of the following holds:
  \begin{enumerate}
    \item There is a \Borel reduction of $E$ to a ccc-idealistic \Borel
      equivalence relation on a \Polish space.
    \item There is a continuous embedding of $\Eone$ into $E$.
  \end{enumerate}  
\end{introtheorem}

As the proof of the \Feldman-\Moore theorem ensures that countable \Borel
equivalence relations on \Polish spaces are countable unions of finite
\Borel subequivalence relations, it is not difficult to see that a \Borel
equivalence relation on a \Polish space is $\sigma$-smooth if and only if it
is a countable union of subequivalence relations that are \Borel reducible to
countable \Borel equivalence relations on \Polish spaces, so Theorem \ref
{introduction:main} also yields:

\begin{introtheorem} \label{introduction:unions}
  Suppose that $E$ is an equivalence relation on a \Polish space that is
  \Borel reducible to a ccc-idealistic \Borel equivalence relation on a \Polish
  space. If $E$ is a countable union of subequivalence relations that are
  \Borel reducible to countable \Borel equivalence relations on \Polish
  spaces, then $E$ is \Borel reducible to a countable \Borel equivalence
  relation on a \Polish space.
\end{introtheorem}

An equivalence relation $E$ on a set $X$ has \definedterm{countable
index} over a subequivalence relation $F$ if every $E$-class is a countable
union of $F$-classes. Generalizing Theorem \ref{introduction:unions} in the
spirit of \cite[Theorem 1.1]{Kittrell}, we show:

\begin{introtheorem} \label{introduction:unions:two}
  Suppose that $E$ is an equivalence relation on a \Polish space that is
  \Borel reducible to a ccc-idealistic \Borel equivalence relation on a \Polish
  space and $\calF$ is a class of strongly-idealistic potentially-\Fsigma
  equivalence relations on \Polish spaces that is closed under countable
  disjoint unions and countable-index \Borel superequivalence relations. If
  $E$ is a countable union of subequivalence relations that are \Borel
  reducible to relations in $\calF$, then $E$ is \Borel reducible to a relation
  in $\calF$.\\
\end{introtheorem}

The \definedterm{saturation} of a set $Y \subseteq X$ with respect to an
equivalence relation $E$ on $X$ is given by $\saturation{Y}{E} \defeq \set
{x \in X}[\exists y \in Y \ x \mathrel{E} y]$, the \definedterm{diagonal} on
$X$ is given by $\diagonal{X} \defeq \set{\pair{x}{y} \in X \times X}[x = y]$,
the product of equivalence relations $E$ and $F$ on $X$ and $Y$ is the
equivalence relation on $X \times Y$ given by $\pair{x}{y} \mathrel{E \times
F} \pair{x'}{y'} \iff (x \mathrel{E} x' \mathand y \mathrel{F} y')$, the
\definedterm{\Friedman--\Stanley jump} of $E$ is the equivalence relation
on $\functions{\N}{X}$ given by $x \mathrel{\jump{E}} y \iff \saturation
{\image{x}{\N}}{E} = \saturation{\image{y}{\N}}{E}$, and we use
$\intersectingrelation{E}$ to denote the binary relation on $\functions{\N}
{X}$ given by $x \mathrel{\intersectingrelation{E}} y \iff \saturation{\image
{x}{\N}}{E} \intersection \saturation{\image{y}{\N}}{E} \neq \emptyset$.
Theorem \ref{introduction:unions:two} yields:

\begin{introtheorem} \label{introduction:intersecting}
  Suppose that $E$ is an equivalence relation on a \Polish space that
  is \Borel reducible to a ccc-idealistic \Borel equivalence relation on a
  \Polish space and $F$ is a strongly-idealistic potentially-\Fsigma
  equivalence relation on a \Polish space. Then the following are equivalent:
  \begin{enumerate}
    \item The equivalence relation $E$ is a countable union of subequivalence
      relations that are \Borel reducible to $F \times \diagonal{\N}$.
    \item There is a \Borel reduction of $E$ to $\intersectingrelation{(F
      \times \diagonal{\N})}$.
    \item There is a \Borel reduction of $E$ to a countable-index \Borel
      superequivalence relation of $F \times \diagonal{\N}$.
    \item There is a \Borel homomorphism from $\pair{E}{\setcomplement
      {E}}$ to $\brokenpair{\jump{(F \times \diagonal{\N})}}{\setcomplement
      {\intersectingrelation{(F \times \diagonal{\N})}}}$.
  \end{enumerate}
  In particular, if there is a \Borel reduction of $E$ to $\intersectingrelation{(F
  \times \diagonal{\N})}$, then there is a \Borel reduction of $E$ to $\jump{(F
  \times \diagonal{\N})}$.
\end{introtheorem}

In \S\ref{secNot}, we introduce basic notation and definitions. In \S\ref
{secCompactOpen}, we review the compact-open topology. In \S\ref
{secCtblIndex}, we establish several preliminaries concerning \Borel
equivalence relations. In \S\ref{secCores}, we characterize the
existence of small definable cores for appropriate families of finite sets.
In \S\ref{secAligned}, we consider approximations to embeddings of
$\Eone$ into itself. In \S\ref{secDichoto}, we establish a pair of technical
dichotomies, and provide an alternate proof of the \Kechris--\Louveau
dichotomy. And in \S\ref{finalSection}, we derive Theorems \ref
{introduction:main}--\ref{introduction:intersecting} from a single common
generalization.

\section{Notation and definitions} \label{secNot}

Given a set $Z$ and a function $f \from X \times Y \to Z$, define
$\horizontalsection{f}{y} \from X \to Z$ and $\verticalsection{f}{x} \from Y \to
Z$ by $\horizontalsection{f}{y}(x) \defeq \verticalsection{f}{x}(y) \defeq f(x,
y)$ for all $x \in X$ and $y \in Y$.

The \definedterm{restriction} of a binary relation $R$ on a set $X$ to a set
$Y \subseteq X$ is the binary relation on $Y$ given by $R \restriction Y
\defeq R \intersection (Y \times Y)$. When considering properties of $R
\restriction Y$ that depend on the ambient space, it will be understood
that this ambient space is $Y$. For instance, if $X$ is a topological
space, then we will say that $R \restriction Y$ is meager, or $R$ is meager
\definedterm{on} $Y$, if it is meager when viewed as a subset of $Y \times
Y$.

Given a topological space $X$ and a property $\phi$ of elements of $X$,
we use $\forcomeagerlymany x \in X \ \phi(x)$ to indicate that $\set{x \in X}
[\phi(x)]$ is comeager, and $\existnonmeagerlymany x \in X \ \phi(x)$ to
indicate that $\set{x \in X}[\phi(x)]$ is not meager.

The \definedterm{complete equivalence relation} on $X$ is given by
$\completeer{X} = X \times X$.

We refer the reader to \cite{Kechris} for basic descriptive set-theoretic
background on \Baire category and analytic sets.

\section{The compact-open topology} \label{secCompactOpen}

Let $X$ and $Y$ be topological spaces. We use $\continuousfunctions{X}
{Y}$ to denote the set of all continuous mappings $f \from X \to Y$. For $A
\subseteq X$ and $B \subseteq Y$, let $\restrictedcontinuousfunctions{A}
{B} \defeq \set{f \in \continuousfunctions{X}{Y}}[\image{f}{A} \subseteq B]$
(we will denote this set by $\restrictedcontinuousfunctions{A}{B}[X][Y]$ in
case of potential ambiguity). We endow $\continuousfunctions{X}{Y}$ with
the \definedterm{compact-open topology}, that is, the topology generated
by the the sets of the form $\restrictedcontinuousfunctions{K}{U}$, where
$K$ ranges over compact subsets of $X$ and $U$ ranges over open
subsets of $Y$.

We summarize below some classical properties of the compact-open
topology. Most of them can be found, for example, in \cite{Engelking}.

\begin{proposition} \label{SubspaceRangeCompactOpen}
  Let $X$ and $Z$ be topological spaces and $Y \subseteq Z$ and view
  $\continuousfunctions{X}{Y}$ as a subset of $\continuousfunctions{X}{Z}$.
  Then the compact-open topology on $\continuousfunctions{X}{Y}$
  coincides with the topology induced by the compact-open topology on
  $\continuousfunctions{X}{Z}$.
\end{proposition}

\begin{propositionproof}
  This immediately follows from the fact that, for $K \subseteq X$ compact
  and $U \subseteq Z$ open, we have $\restrictedcontinuousfunctions{K}{U
  \intersection Y}[X][Y] = \restrictedcontinuousfunctions{K}{U}[X][Z]
  \intersection \continuousfunctions{X}{Y}$.
\end{propositionproof}

\begin{proposition}[{see \cite[Theorem 3.4.2]{Engelking}}]
  \label{ContCompo}
  Let $X$, $Y$, and $Z$ be topological spaces, with $Y$ being locally
  compact. Then the composition mapping $\continuousfunctions{Y}{Z}
  \times \continuousfunctions{X}{Y} \to \continuousfunctions{X}{Z}$, given
  by $\pair{f}{g} \mapsto f \composition g$, is continuous.
\end{proposition}

\begin{proposition} \label{ContEval}
  Let $Y$ and $Z$ be topological spaces, with $Y$ being locally compact.
  Then the evaluation mapping $\continuousfunctions{Y}{Z} \times Y \to Z$,
  given by $\pair{f}{x} \mapsto f(x)$, is continuous.
\end{proposition}

\begin{propositionproof}
  Apply Proposition \ref{ContCompo} with $X$ being a singleton.
\end{propositionproof}

\begin{proposition}[{see \cite[Theorems 3.4.3 and 3.4.8]{Engelking}}]
  \label{ExpCompactOpen}
  Let $X$, $Y$, and $Z$ be topological spaces, with $X$ being locally
  compact and $Y$ being Hausdorff.
  \begin{enumerate}
    \item For every $f \in \continuousfunctions{X \times Y}{Z}$, the mapping
      $\Lambda(f) \from Y \to \continuousfunctions{X}{Z}$, given by $y
        \mapsto \horizontalsection{f}{y}$, is continuous.
    \item The mapping $\Lambda \from \continuousfunctions{X \times Y}{Z}
      \to \continuousfunctions{Y}{\continuousfunctions{X}{Z}}$ hence defined
        is a homeomorphism.
  \end{enumerate}
\end{proposition}

\begin{proposition}[{see \cite[Exercise 3.4.B]{Engelking}}]
  \label{ProductCompactOpen}
  Let $X_1$, $X_2$, $Y_1$, and $Y_2$ be topological spaces, with $X_1$
  and $X_2$ being Hausdorff. Then the mapping $\continuousfunctions{X_1}
  {Y_1} \times \continuousfunctions{X_2}{Y_2} \to \continuousfunctions{X_1
  \times X_2}{Y_1 \times Y_2}$, given by $\pair{f_1}{f_2} \mapsto f_1 \times
  f_2$, is a homeomorphic embedding.
\end{proposition}

Let $X$ be a compact topological space and $\pair{Y}{d}$ be a metric
space. The \definedterm{uniform metric} on $\continuousfunctions{X}{Y}$
associated with $d$ is given by $\uniformmetric{d}(f, g) \defeq \sup_{x \in
X} d(f(x), g(x))$.

\begin{proposition}[{see \cite[Theorems 4.2.17 and 4.3.13]{Engelking}}]
  \label{PropUnifMetric}
  Let $X$ be a compact topological space and $\pair{Y}{d}$ be a metric
  space. Then $\uniformmetric{d}$ is a metric on $\continuousfunctions{X}
  {Y}$ that is compatible with the compact-open topology. Moreover, if
  $\pair{Y}{d}$ is complete, then so too is $\pair{\continuousfunctions{X}
  {Y}}{\uniformmetric{d}}$.
\end{proposition}

\begin{proposition}[{see \cite[Theorem 3.4.16 and Exercise 4.3.F]
  {Engelking}}] \label{PolishCompactOpen}
  Let $X$ and $Y$ be topological spaces, with $X$ being locally compact
  and second countable.
  \begin{enumerate}
    \item If $Y$ is second countable, then $\continuousfunctions{X}{Y}$ is
      second countable.
    \item If $Y$ is completely metrizable, then $\continuousfunctions{X}{Y}$
      is completely metrizable.
  \end{enumerate}
  In particular, if $Y$ is \Polish, then $\continuousfunctions{X}{Y}$ is \Polish.
\end{proposition}

Given topological spaces $X$, $Y$, and $Z$ and $\calA \subseteq
\continuousfunctions{X}{Z}$, we denote by $\sectionpower{\calA}{Y}$ the
set of all $f \in \continuousfunctions{X \times Y}{Z}$ such that, for all $y \in
Y$, $\horizontalsection{f}{y} \in \calA$.

\begin{lemma} \label{ArrowPreservesOpen}
  Let $X$, $Y$, and $Z$ be topological spaces, with $X$ being locally
  compact and $Y$ being compact. Let $\calU \subseteq
  \continuousfunctions{X}{Z}$ be an open set. Then $\sectionpower{\calU}
  {Y}$ is an open subset of $\continuousfunctions{X \times Y}{Z}$.
\end{lemma}

\begin{lemmaproof}
  Keeping the notation from the statement of Proposition \ref
  {ExpCompactOpen}, we have $\sectionpower{\calU}{Y} = \preimage
  {\Lambda}{\restrictedcontinuousfunctions{Y}{\calU}}$, hence Proposition
  \ref{ExpCompactOpen} immediately gives the desired result.
\end{lemmaproof}

Given a natural number $n$, a set $X$, and a sequence $\sequence{Y_i}[i
\in \N]$, we identify $(X \times \product[i < n][Y_i]) \times \product[i \ge n]
[Y_i]$ with $X \times \product[i \in \N][Y_i]$.

\begin{lemma} \label{CompatibleMetrics}
  Let $X$ be a locally-compact second-countable topological space,
  $\sequence{Y_n}[n \in \N]$ be a sequence of compact topological spaces,
  and $Z$ be a completely metrizable topological space. Then there are
  compatible complete metrics $d_n$ on $\continuousfunctions{X \times
  \product[i < n][Y_i]}{Z}$ such that $\diam_\N(\sectionpower{\calA}{\product
  [i \ge n][Y_i]}) \le \diam_n(\calA)$ for all $n \in \N$ and $\calA \subseteq
  \continuousfunctions{X \times \product[i < n][Y_i]}{Z}$, where $\diam_n$
  denotes the diameter relative to $d_n$ for all $n \in \N \union \set{\N}$.
\end{lemma}

\begin{lemmaproof}
  We first deal with the special case when $X$ is a singleton. Fix a
  compatible complete metric $d$ on $Z$. For all $n \in \N \union \set{\N}$,
  identify $X \times \product[i < n][Y_i]$ with $\product[i < n][Y_i]$ and let
  $d_n$ be the uniform metric on $\continuousfunctions{\product[i < n][Y_i]}
  {Z}$. Suppose now that $n \in \N$ and $\calA \subseteq
  \continuousfunctions{\product[i < n][Y_i]}{Z}$, and observe that if $f, g \in
  \sectionpower{\calA}{\product[i \ge n][Y_i]}$, then
  \begin{align*}
    d_\N(f, g)
      & = \sup_{y \in \product[i \in \N][Y_i]} d(f(y), g(y)) \\
      & = \sup_{v \in \product[i \ge n][Y_i]} \sup_{u \in \product
        [i < n][Y_i]} d(\horizontalsection{f}{v}(u), \horizontalsection{g}{v}(u)) \\
      & = \sup_{v \in \product[i \ge n][Y_i]} d_n(\horizontalsection{f}{v},
        \horizontalsection{g}{v}) \\
      & \le \diam_n(\calA),
  \end{align*}
  so $\diam_\N(\sectionpower{\calA}{\product[i \ge n][Y_i]}) \le \diam_n
  (\calA)$.

  We now deal with the general case. Let $Z' = \continuousfunctions{X}{Z}$.
  By Proposition \ref{ExpCompactOpen}, for every $n \in \N \union
  \set{\N}$, the mapping $\Lambda_n \from  \continuousfunctions{X \times
  \product[i < n][Y_i]}{Z} \to \continuousfunctions{\product[i < n][Y_i]}{Z'}$,
  defined by $\Lambda_n(f)(y) \defeq \horizontalsection{f}{y}$, is a
  homeomorphism. Moreover, it is easy to verify that, for every $f \in
  \continuousfunctions{X \times \product[i \in \N][Y_i]}{Z}$, $n \in \N$, and $y
  \in \product[i \ge n][Y_i]$, we have $\Lambda_n(\horizontalsection
  {f}{y}) = \horizontalsection{\Lambda_\N(f)}{y}$. It follows that, for every
  $n \in \N$ and $\calA \subseteq \continuousfunctions{X \times \product[i < n]
  [Y_i]}{Z}$, we have $\image{\Lambda_\N}{\sectionpower{\calA}{\product[i \ge n]
  [Y_i]}} = \sectionpower{\image{\Lambda_n}{\calA}}{\product[i \ge n][Y_i]}$.

  By Proposition \ref{PolishCompactOpen}, $Z'$ is completely metrizable, so
  we can apply the special case to find, for every $n \in \N
  \union \set{\N}$, a complete metric $d_n'$ on  $\continuousfunctions
  {\product[i < n][Y_i]}{Z'}$ such that, for every $n \in \N$ and $\calA'
  \subseteq \continuousfunctions{\product[i < n][Y_i]}{Z'}$, we have
  $\diam_\N(\sectionpower{(\calA')}{\product[i \ge n][Y_i]}) \le \diam_n(\calA')$.
  We define, for every $n \in \N \union \set{\N}$, the metric $d_n$ on
  $\continuousfunctions{X \times \product[i < n][Y_i]}{Z}$ as the pullback of
  $d_n'$ through $\Lambda_n$. Hence, for $n \in \N$ and $\calA \subseteq
  \continuousfunctions{X \times \product[i < n][Y_i]}{Z}$, we have
  \begin{align*}
    \diam_\N(\sectionpower{\calA}{\product[i \ge n][Y_i]})
      & = \diam_\N(\image{\Lambda_\N}{\sectionpower{\calA}{\product[i \ge n]
        [Y_i]}}) \\
      & = \diam_\N(\sectionpower{\image{\Lambda_n}{\calA}}{\product[i \ge n]
        [Y_i]}) \\
      & \le \diam_n(\image{\Lambda_n}{\calA}) \\
      & = \diam_n(\calA),
  \end{align*}
  which completes the proof.
\end{lemmaproof}

\begin{lemma} \label{BasisCompactOpen}
  Let $X$ be a zero-dimensional compact \Polish space, $\Omega \subseteq
  X$ be an open subset, $Y$ be a topological space, and $\calB$ be a
  basis of open subsets of $Y$ that is closed under finite intersection. Then
  the sets of the form $\intersection[i < n][\restrictedcontinuousfunctions{K_i}
  {U_i}]$, where $n \in \N$, $\sequence{K_i}[i < n]$ is a sequence of
  nonempty pairwise disjoint clopen subsets of $X$ that are contained in
  $\Omega$, and $\sequence{U_i}[i < n]$ is a sequence of nonempty
  elements of $\calB$, form a basis of nonempty open subsets of
  $\continuousfunctions{\Omega}{Y}$.
  
  In the special case when $\Omega = X$, the sets as above, but where we,
  moreover, require that $\sequence{K_i}[i < n]$ is a partition of $X$, form a
  basis of nonempty open subsets of $\continuousfunctions{X}{Y}$.
\end{lemma}

\begin{lemmaproof}
  Keeping the notation from the statement of the lemma and taking $y_i \in
  U_i$ for all $i < n$ when $n > 0$, the mapping $f \from \Omega \to Y$,
  defined by $f(x) = y_i$ for all $x \in K_i$ and $f(x) = y_0$ for all $x \notin
  \union[i < n][K_i]$, is an element of $\intersection[i < n]
  [\restrictedcontinuousfunctions{K_i}{U_i}]$, which is hence nonempty. We
  now show that these sets form a basis of open sets of
  $\continuousfunctions{\Omega}{Y}$.

  Let $\calU$ be an open subset of $\continuousfunctions{\Omega}{Y}$ and
  $f \in \calU$. We can find $k \in \N$, nonempty compact sets $L_0,
  \ldots, L_{k-1} \subseteq \Omega$, and open sets $V_0, \ldots, V_{k-1}
  \subseteq Y$ such that $f \in \intersection[i < k]
  [\restrictedcontinuousfunctions{L_i}{V_i}] \subseteq \calU$. In particular, for
  every $i < k$, we have $\image{f}{L_i} \subseteq V_i$. Fix $i < k$. For
  every $x \in L_i$, there exists $W_i^x \in \calB$ such that $f(x) \in W_i^x
  \subseteq V_i$. We can find a clopen subset $C_i^x$ of $X$, contained in
  $\Omega$, such that $x \in C_i^x \subseteq \preimage{f}{W_i^x}$. Since
  $L_i$ is compact, we can find $x_i^0, \ldots, x_i^{l_i} \in L_i$ such that
  $L_i \subseteq \union[j \le l_i][C_i^{x_i^j}]$. Then
  \begin{equation*}
    f \in \intersection[i < k][{\intersection[j \le l_i][\restrictedcontinuousfunctions
      {C_i^{x_i^j}}{W_i^{x_i^j}}]}] \subseteq \calU.
  \end{equation*}
  This shows that, by shrinking $\calU$ if necessary, we can assume that the
  $L_i$'s are clopen in $X$ and the $V_i$'s are elements of $\calB$.

  For every set $s \subseteq k$, let $K_s = (\intersection[i \in s][L_i])
  \setminus (\union[i \in k \setminus s][L_i])$ and $U_s = \intersection[i \in
  s][V_i]$. Let $S = \set{s \subseteq k}[K_s \neq \emptyset] \setminus \set
  {\emptyset}$. For every $s \in S$, $K_s$ is contained in one of the
  $L_i$'s, and is therefore a subset of $\Omega$; clearly, the $K_s$'s are
  pairwise disjoint and clopen. Moreover, we have $\intersection[s \in S]
  [\restrictedcontinuousfunctions{K_s}{U_s}] = \intersection[i < k]
  [\restrictedcontinuousfunctions{L_i}{V_i}]$. Thus the latter set is an open
  neighborhood of $f$ of the desired form that is contained in $\calU$.

  In the special case when $\Omega = X$, a basic open subset $\union[i <
  n][\restrictedcontinuousfunctions{K_i}{U_i}]$ of $\continuousfunctions{X}
  {Y}$, where $n \in \N$, $\sequence{K_i}[i < n]$ is a sequence of nonempty
  pairwise disjoint clopen subsets of $X$, and $\sequence{U_i}[i < n]$ is a
  sequence of nonempty elements of $\calB$, can be rewritten as $\union[i
  \le n][\restrictedcontinuousfunctions{K_i}{U_i}]$, where $K_n = X
  \setminus \union[i < n][K_i]$ and $U_n = X$.
\end{lemmaproof}

Let $X$ and $Y$ be topological spaces and $\calU$ be an open subset of
$\continuousfunctions{X}{Y}$. We say that $\calU$ is \definedterm{right
stable} when there exists a neighborhood $\calV$ of the identity in
$\continuousfunctions{X}{X}$ such that $\calU \composition \calV = \calU$.

\begin{corollary} \label{BasisRightStable}
  Let $X$ be a zero-dimensional second-countable compact topological
  space, $\Omega \subseteq X$ be an open set, and $Y$ be a
  second-countable topological space. Then $\continuousfunctions
  {\Omega}{Y}$ admits a countable basis of open subsets whose
  elements are right stable.
\end{corollary}

\begin{corollaryproof}
  Fix a countable basis $\calB$ for $Y$ that is closed under finite
  intersection. Then, by Lemma \ref{BasisCompactOpen}, the sets of the
  form $\calU \defeq \intersection[i < n][\restrictedcontinuousfunctions{K_i}
  {U_i}]$, where $n \in \N$, $\sequence{K_i}[i < n]$ is a sequence of
  nonempty pairwise disjoint clopen subsets of $X$ that are contained in
  $\Omega$, and $\sequence{U_i}[i < n]$ is a sequence of nonempty
  elements of $\calB$, form a basis for $\continuousfunctions{\Omega}{Y}$.
  Moreover, this basis is countable. It remains to show that every such set
  is right stable. Keeping the notation as above, the set $\calV \defeq
  \intersection[i < n][\restrictedcontinuousfunctions{K_i}{K_i}]$ is an open
  neighborhood of the identity in $\continuousfunctions{\Omega}{\Omega}$
  and $\calU \composition \calV = \calU$.
\end{corollaryproof}

\section{Preliminaries} \label{secCtblIndex}

Clearly the family of sets $Y \subseteq X$ on which $E$ has countable
index over $F$ is closed under containment and $F$-saturation.

\begin{lemma} \label{PushForwardCountableIndex}
  Suppose that $X$ and $Y$ are sets, $E \subseteq E'$ are equivalence
  relations on $X$, $F \subseteq F'$ are equivalence relations on $Y$, and
  $\phi \from X \to Y$ is a homomorphism from $\pair{E}{\setcomplement
  {E'}}$ to $\pair{F}{\setcomplement{F'}}$. If $E'$ has countable index over
  $E$, then $F'$ has countable index over $F$ on $\image{\phi}{X}$.
\end{lemma}

\begin{lemmaproof}
  It is sufficient to show that if $\phi(x) = y$ and $\equivalenceclass{x}{E'} =
  \union[n \in \N][\equivalenceclass{x_n}{E}]$, then $\equivalenceclass{y}{F'
  \restriction \image{\phi}{X}} \subseteq \union[n \in \N][\equivalenceclass
  {\phi(x_n)}{F}]$. Towards this end, suppose that $y' \in \equivalenceclass
  {y}{F' \restriction \image{\phi}{X}}$, and fix $x' \in X$ for which $\phi(x') =
  y'$. As $\phi$ is a cohomomorphism from $E'$ to $F'$, it follows that $x
  \mathrel{E'} x'$, so there exists $n \in \N$ for which $x_n \mathrel{E} x'$,
  and since $\phi$ is a homomorphism from $E$ to $F$, it follows that $y'
  \in \equivalenceclass{\phi(x_n)}{F}$.
\end{lemmaproof}

The following fact will prove useful in complexity calculations:

\begin{proposition} \label{CalculationComplexityCountableIndex}
  Suppose that $X$ is a \Polish space, $E$ is
  an analytic equivalence relation on $X$, and $F$ is a co-analytic
  equivalence relation on $X$. Then exactly one of the following holds:
  \begin{enumerate}
    \item The equivalence relation $E$ has countable index over $E
      \intersection F$.
    \item There is a continuous function $\phi \from \Cantorspace \to X$
      with the property that $\preimage{(\phi \times \phi)}{E \setminus F}$ is
      comeager.
  \end{enumerate}
\end{proposition}

\begin{propositionproof}
  To see $(2) \implies \neg(1)$, suppose that $R \defeq \preimage{(\phi
  \times \phi)}{E \setminus F}$ is comeager, appeal to \Mycielski's theorem
  (see, for example, \cite[Theorem 19.1]{Kechris}) to obtain a continuous
  homomorphism $\psi \from \Cantorspace \to \Cantorspace$ from
  $\setcomplement{\diagonal{\Cantorspace}}$ to $R$, and observe that
  $\phi \composition \psi$ is a continuous homomorphism from
  $\setcomplement{\diagonal{\Cantorspace}}$ to $E \setminus F$, thus
  condition (1) fails.

  To see $\neg(1) \implies (2)$, suppose that there exists $x \in X$ for which
  $F \restriction \equivalenceclass{x}{E}$ has uncountably many classes,
  appeal to the straightforward generalization of \Silver's perfect set
  theorem (see \cite{SilverEqRel}) to analytic subsets of \Polish spaces to
  obtain a continuous embedding $\phi \from \Cantorspace \to
  \equivalenceclass{x}{E}$ of $\diagonal{\Cantorspace}$ into $F \restriction
  \equivalenceclass{x}{E}$, and observe that $\preimage{(\phi \times \phi)}
  {E \setminus F} = \setcomplement{\diagonal{\Cantorspace}}$, thus
  condition (2) holds.
\end{propositionproof}

Given a binary relation $R$ on a set $X$, let $\generateder{R}$ denote the
smallest equivalence relation on $X$ containing $R$.

\begin{proposition} \label{ComplexityGeneratedCountableIndex}
  Suppose that $X$ and $Y$ are \Polish spaces, $R \subseteq (X \times X)
  \times Y$ is analytic, and $F$ is a co-analytic equivalence relation on $X$.
  Then $\set{y \in Y}[\generateder{\horizontalsection{R}{y}} \text{ has
  countable index over } \generateder{\horizontalsection{R}{y}} \intersection
  F]$ is co-analytic.
\end{proposition}

\begin{propositionproof}
  If $y \in Y$, then Proposition \ref{CalculationComplexityCountableIndex}
  ensures that $\generateder{\horizontalsection{R}{y}}$ has countable index
  over $\generateder{\horizontalsection{R}{y}} \intersection F$ if and only if
  there is no continuous function $\phi \from \Cantorspace \to X$ for which
  $\preimage{(\phi \times \phi)}{\generateder{\horizontalsection{R}{y}}
  \setminus F}$ is comeager. But Proposition \ref{ContEval} and \cite
  [Theorem 29.22]{Kechris} imply that the set of $y$ for which this holds is
  co-analytic.
\end{propositionproof}

\begin{corollary} \label{ComplexityCountableIndex}
  Suppose that $X$ and $Y$ are \Polish spaces, $E$ is an analytic
  equivalence relation on $X$, $F$ is a co-analytic equivalence relation on
  $X$, and $R \subseteq X \times Y$ is an analytic set. Then $\set{y \in Y}
  [E \text{ has countable index over } E \intersection F$ on
  $\horizontalsection{R}{y}]$ is co-analytic.
\end{corollary}

\begin{corollaryproof}
  Set $S \defeq \set{\pair{\pair{w}{x}}{y} \in E \times Y}[w, x \in
  \horizontalsection{R}{y}]$. If $y \in Y$, then $\generateder
  {\horizontalsection{S}{y}} = (E \restriction \horizontalsection{R}{y}) \union
  \diagonal{X}$, so $\generateder{\horizontalsection{S}{y}}$ has countable
  index over $\generateder{\horizontalsection{S}{y}} \intersection F$ if and
  only if $E$ has countable index over
  $E \intersection F$ on $\horizontalsection{R}{y}$, thus Proposition 
  \ref{ComplexityGeneratedCountableIndex} yields the desired result.
\end{corollaryproof}

\begin{corollary} \label{ComplexityCountablyManyClasses}
  Suppose that $X$ and $Y$ are \Polish spaces, $F$ is a co-analytic
  equivalence relation on $X$, and $R \subseteq X \times Y$ is analytic.
  Then $\set{y \in Y}[F \text{ has only countably many classes on }
  \horizontalsection{R}{y}]$ is co-analytic.
\end{corollary}

\begin{corollaryproof}
  If $y \in Y$, then $\completeer{X}$ has countable index over $F$ on
  $\horizontalsection{R}{y}$ if and only if $F$ has only countably many
  classes on $\horizontalsection{R}{y}$, so Corollary \ref
  {ComplexityCountableIndex} yields the desired result.
\end{corollaryproof}

We say that a property $\Phi$ of subsets of a \Polish space $X$ is
\definedterm{$\Piclass[1][1]$-on-$\Sigmaclass[1][1]$} if $\set{x \in X}[\Phi
(\horizontalsection{R}{y})]$ is co-analytic for every \Polish space $Y$ and
analytic set $R \subseteq X \times Y$.

\begin{corollary} \label{ReflexionCountableIndexPartial}
  Suppose that $X$ is a \Polish space, $E$ is an analytic equivalence
  relation on $X$, $F$ is a co-analytic equivalence relation on $X$, and
  $A \subseteq X$ is an analytic set on which $E$ has countable index
  over $E \intersection F$. Then there is a \Borel set $B \supseteq A$ on
  which $E$ has countable index over $E \intersection F$.
\end{corollary}

\begin{corollaryproof}
  By Corollary \ref{ComplexityCountableIndex}, the property (of $A$) that
  $E \restriction A$ has countable index over $(E \intersection F) \restriction
  A$ is $\Piclass[1][1]$-on-$\Sigmaclass[1][1]$, so the first reflection
  theorem (see, for example, \cite[Theorem 35.10]{Kechris}) yields the
  desired result.
\end{corollaryproof}

\begin{proposition} \label{ReflexionCountableIndex}
  Suppose that $X$ is a \Polish space, $E$ is an analytic equivalence
  relation on $X$, $F$ is a \Borel equivalence relation on $X$ contained in
  $E$, and $A \subseteq X$ is an analytic set on which $E$ has countable
  index over $F$. Then there is an $F$-invariant \Borel set $B \supseteq A$
  on which $E$ has countable index over $F$.
\end{proposition}

\begin{propositionproof}
  Set $A_0 \defeq \saturation{A}{F}$. Given an analytic set $A_n \subseteq
  X$ on which $E$ has countable index over $F$, appeal to Corollary \ref
  {ReflexionCountableIndexPartial} to obtain a \Borel set $B_n \supseteq
  A_n$ on which $E$ has countable index over $F$, and set $A_{n+1}
  \defeq \saturation{B_n}{F}$. Define $B \defeq \union[n \in \N][A_n] =
  \union[n \in \N][B_n]$.
\end{propositionproof}

\begin{proposition} \label{SwitchSidesCountableIndex}
  Suppose that $X$ and $Y$ are \Polish spaces, $E$ and $F$ are
  \Borel equivalence relations on $X$ and $Y$, $\phi \from X \to Y$ is
  a \Borel reduction of $E$ to $F$, and $E'$ is a countable-index \Borel
  superequivalence relation of $E$. Then there is a countable-index
  \Borel superequivalence relation $F'$ of $F$ for which $\phi$ is a
  reduction of $E'$ to $F'$.
\end{proposition}

\begin{propositionproof}
  Define $F_0 \defeq \image{(\phi \times \phi)}{E'}$.
  
  \begin{claim}
    The relation $F_0$ is transitive.
  \end{claim} 
  
  \begin{claimproof}
    Suppose that $y_1 \mathrel{F_0} y_2 \mathrel{F_0} y_3$. Then there
    exist $\pair{w_1}{w_2}, \pair{x_2}{x_3} \in E'$ such that $\phi(w_i) =
    y_i$ for all $i \in \set{1, 2}$ and $\phi(x_i) = y_i$ for all $i \in \set{2, 3}$.
    As $\phi$ is a cohomomorphism from $E$ to $F$, it follows that
    $w_1 \mathrel{E'} w_2 \mathrel{E} x_2 \mathrel{E'} x_3$, so $w_1
    \mathrel{E'} x_3$, thus $y_1 \mathrel{F_0} y_3$.
  \end{claimproof}
    
  It follows that $F_0$ is an equivalence relation on $\image{\phi}{X}$.
  
  \begin{claim}
    The function $\phi$ is a cohomomorphism from $E'$ to $F_0$.
  \end{claim}
  
  \begin{claimproof}
    Suppose that $\phi(x_1) \mathrel{F_0} \phi(x_2)$. Then there exists
    $\pair{w_1}{w_2} \in E'$ such that $\phi(w_i) = \phi(x_i)$ for all $i \in \set
    {1, 2}$. As $\phi$ is a cohomomorphism from $E$ to $F$, it follows that
    $x_1 \mathrel{E} w_1 \mathrel{E'} w_2 \mathrel{E} x_2$, so $x_1
    \mathrel{E'} x_2$.
  \end{claimproof}
  
  Set $F_0' \defeq \generateder{F \union F_0}$, and note that $F_0 = F_0'
  \restriction \image{\phi}{X}$, so $\phi$ is a cohomomorphism from $E'$ to
  $F_0'$. As $E'$ has countable index over $E$, Proposition \ref
  {PushForwardCountableIndex} ensures that $F_0$ has countable index
  over $F \restriction \image{\phi}{X}$, and since $\equivalenceclass{y}{F}
  \neq \equivalenceclass{y}{F_0'} \implies \equivalenceclass{y}{F}
  \intersection \image{\phi}{X} \neq \emptyset$ for all $y \in Y$, it follows
  that $F_0'$ has countable index over $F$.
  
  Suppose now that $n \in \N$ and $F_n'$ is an analytic superequivalence
  relation of $F_0'$ such that $\phi$ is a cohomomorphism from $E'$ to
  $\generateder{F_n'}$ and $\generateder{F_n'}$ has countable index over
  $F$. As Proposition \ref{ComplexityGeneratedCountableIndex} ensures
  that this is a $\Piclass[1][1]$-on-$\Sigmaclass[1][1]$ property (of $F_n'$),
  the first reflection theorem yields a \Borel superrelation $R_n$ of $F_n'$
  for which $\phi$ is a cohomomorphism from $E'$ to $\generateder{R_n}$
  and $\generateder{R_n}$ has countable index over $F$. Set $F_{n+1}'
  \defeq \generateder{R_n}$.
  
  Clearly $F' \defeq \union[n \in \N][F_n'] = \union[n \in \N][R_n]$ is a \Borel
  superequivalence relation of $F$. To see that $\phi$ is a
  cohomomorphism from $E'$ to $F'$, and therefore a reduction of $E'$ to
  $F'$, note that if $\phi(w) \mathrel{F'} \phi(x)$, then there exists $n \in \N$
  for which $\phi(w) \mathrel{F_n'} \phi(x)$, so $w \mathrel{E'} x$. To see
  that $F'$ has countable index over $F$, note that $\equivalenceclass{y}
  {F'} = \union[n \in \N][\equivalenceclass{y}{F_n'}]$ for all $y \in Y$ and $F$
  has only countably many classes on $\equivalenceclass{y}{F_n'}$ for all
  $n \in \N$ and $y \in Y$.
\end{propositionproof}

For equivalence relations $E$ and $F$ on sets $X$ and $Y$, we say that a
set $R \subseteq X \times Y$ \definedterm{induces a partial injection} of
$X / E$ into $Y / F$ if $x \mathrel{E} x' \iff y \mathrel{F} y'$ for all $\pair{x}
{y}, \pair{x'}{y'} \in R$.

\begin{proposition} \label{InducesAnInjection}
  Suppose that $X$ and $Y$ are \Polish spaces, $E$ and $F$ are \Borel
  equivalence relations on $X$ and $Y$, and $R \subseteq X \times Y$ is
  an analytic set inducing a partial injection of $X / E$ into $Y / F$. Then
  there is an $(E \times F)$-invariant \Borel set $S \supseteq R$ inducing a
  partial injection of $X / E$ into $Y / F$.
\end{proposition}

\begin{propositionproof}
  Set $R_0 \defeq R$. Given an analytic set $R_n \subseteq X \times Y$
  inducing a partial injection of $X / E$ into $Y / F$, appeal to the first
  reflection theorem to obtain a \Borel set $S_n \supseteq R_n$ inducing
  a partial injection of $X / E$ into $Y / F$, and set $R_{n+1} \defeq
  \saturation{S_n}{E \times F}$. Define $S \defeq \union[n \in \N][R_n] =
  \union[n \in \N][S_n]$.
\end{propositionproof}

As a first corollary of this result, we have the following:

\begin{proposition} \label{InvariantReduction}
  Suppose that $X$ and $Y$ are \Polish spaces, $E$ is a \Borel
  equivalence relation on $X$, $F$ is a strongly-idealistic \Borel
  equivalence relation on $Y$, and $A \subseteq X$ is an analytic set on
  which $E$ is \Borel reducible to $F$. Then there is an $E$-invariant \Borel
  set $B \supseteq A$ on which $E$ is \Borel reducible to $F$.
\end{proposition}

\begin{propositionproof}
  As the graph of any \Borel reduction $\pi \from A \to Y$ of $E \restriction
  A$ to $F$ is necessarily analytic (see, for example, the proof of \cite
  [Proposition 12.4]{Kechris}), the proposition follows from an application of
  Proposition \ref{InducesAnInjection} and \cite[Proposition 2.8]
  {deRancourtMiller} to the graph of $\pi$.
\end{propositionproof}

\begin{proposition} \label{CountableIndexToIntersectionReduction}
  Suppose that $X$ and $Y$ are \Polish spaces, $E$ is an equivalence
  relation on $X$, $F$ is a strongly-idealistic \Borel equivalence relation on
  $Y$, and there is a \Borel reduction $\pi \from X \to Y$ of $E$ to a
  countable-index \Borel superequivalence relation $F'$ of $F$. Then there
  is a \Borel homomorphism $\phi \from X \to \functions{\N}{Y}$ from $\pair
  {E}{\setcomplement{E}}$ to $\pair{\jump{F}}{\setcomplement
  {\intersectingrelation{F}}}$.
\end{proposition}

\begin{propositionproof}
  By \cite[Theorem 2.12]{deRancourtMiller}, there are \Borel functions
  $\phi_n \from Y \to Y$ such that $\equivalenceclass{y}{F'} = \union[n \in \N]
  [\equivalenceclass{\phi_n(y)}{F}]$ for all $y \in Y$. Then the function $\phi
  \from X \to \functions{\N}{Y}$, given by $\phi(x)(n) = (\phi_n \composition
  \pi)(x)$, is as desired.
\end{propositionproof}

\begin{remark}
  Proposition \ref{SwitchSidesCountableIndex} easily implies the
  generalization of the converse of Proposition \ref
  {CountableIndexToIntersectionReduction} in which $F$ need not be
  strongly idealistic.
\end{remark}

We will use Proposition \ref{CountableIndexToIntersectionReduction} in
conjunction with the following:

\begin{proposition} \label{IntersectingCountableUnion}
  Suppose that $X$ and $Y$ are \Polish spaces, $E$ and $F$ are \Borel
  equivalence relations on $X$ and $Y$, and there is a \Borel reduction $\pi
  \from X \to \functions{\N}{Y}$ of $E$ to $\intersectingrelation{F}$. Then
  $E$ is a countable union of subequivalence relations that are \Borel
  reducible to $F$.
\end{proposition}

\begin{propositionproof}
  We can assume, without loss of generality, that $X$ is a subset of
  $\Cantorspace$ and that the \Borel subsets of $X$ are the \Borel subsets
  of $\Cantorspace$ contained in $X$. Define $\phi \from X \to \functions{\N
  \times \N \times \N}{Y}$ by
  \begin{equation*}
    \phi(x)(i,j,k) =
      \begin{cases}
        \pi(x)(i) & \text{if $x(k) = 0$ and} \\
        \pi(x)(j) & \text{if $x(k) = 1$.}
      \end{cases}
  \end{equation*}
  For all $i, j, k \in \N$, let $E_{i,j,k}$ be the equivalence relation on $X$
  given by $x \mathrel{E_{i,j,k}} y \iff \phi(x)(i,j,k) \mathrel{F} \phi(y)(i,j,k)$.
  The fact that $\pi$ is a cohomomorphism from $E$ to
  $\intersectingrelation{F}$ ensures that each $E_{i,j,k}$ is a
  subequivalence relation of $E$. To see that the union of these
  equivalence relations is $E$, suppose that $x$ and $y$ are distinct
  $E$-related points of $X$, and fix $k \in \N$ for which $x(k) \neq y(k)$. By
  reversing the roles of $x$ and $y$ if necessary, we can assume that $x(k)
  = 0$. As $\pi$ is a homomorphism from $E$ to $\intersectingrelation{F}$,
  there exist $i,j \in \N$ for which $\pi(x)(i) \mathrel{F} \pi(y)(j)$, in which
  case $\phi(x)(i,j,k) = \pi(x)(i) \mathrel{F} \pi(y)(j) = \phi(y)(i,j,k)$, thus $x
  \mathrel{E_{i,j,k}} y$.
\end{propositionproof}

\begin{remark}
  Conversely, a straightforward argument shows that if $E$ is a countable
  union of subequivalence relations that are \Borel reducible to $F \times
  \diagonal{\N}$, then $E$ is \Borel reducible to $\intersectingrelation{(F
  \times \diagonal{\N})}$.
\end{remark}

We also obtain a useful closure property of the class of strongly-idealistic
potentially-\Fsigma \Borel equivalence relations on \Polish spaces:

\begin{proposition} \label{idealistic:countableindex}
  Suppose that $X$ is a \Polish space and  $E$ is a \Borel equivalence
  relation on $X$ that has a countable-index strongly(-ccc)-idealistic
  potentially-\Fsigma subequivalence relation $F$. Then $E$ is strongly
  (ccc) idealistic and potentially \Fsigma.
\end{proposition}

\begin{propositionproof}
  By \cite[Theorem 2.12]{deRancourtMiller}, there are \Borel functions
  $\phi_n \from X \to X$ with the property that $\equivalenceclass{x}{E} =
  \union[n \in \N][\equivalenceclass{\phi_n(x)}{F}]$ for all $x \in X$.
  
  To see that $E$ is potentially \Fsigma, appeal to standard change of
  topology results to obtain a \Polish topology on $X$, generating the same
  \Borel sets as the given topology, with respect to which $F$ is \Fsigma
  and $\phi_n$ is continuous for all $n \in \N$. As $E = \union[n \in \N]
  [\preimage{(\phi_n \times \id{X})}{F}]$, it is \Fsigma with respect to any
  such topology.
  
  To see that $E$ is strongly (ccc) idealistic, fix a witness $x \mapsto
  \calJ_x$ to the strong (ccc) idealisticity of $F$.
  
  \begin{claim} \label{idealistic:BorelonBorel}
    Suppose that $\phi \from X \to X$ is \Borel. Then $x \mapsto
    \calJ_{\phi(x)}$ is strongly \Borel-on-\Borel.
  \end{claim}

  \begin{claimproof}
    Given a \Polish space $Y$ and a \Borel set $R \subseteq (X \times Y)
    \times X$, set $S \defeq \set{\pair{\pair{x}{\pair{w}{y}}}{z} \in (X \times (X
    \times Y)) \times X}[\phi(w) = x \mathand \pair{w}{y} \mathrel{R} z]$,
    define $\psi \from X \times Y \to X \times (X \times Y)$ by $\psi(w, y)
    \defeq \pair{\phi(w)}{\pair{w}{y}}$, and observe that
    \begin{align*}
      \set{& \pair{w}{y} \in X \times Y}[\verticalsection{R}{\pair{w}{y}} \in
        \calI_{\phi(w)}]  \\
          & = \preimage{\psi}{\set{\pair{\phi(w)}{\pair{w}{y}}}[\pair{w}{y} \in X
            \times Y \mathand \verticalsection{R}{\pair{w}{y}} \in \calI_{\phi(w)}]}
            \\
          & = \preimage{\psi}{\set{\pair{x}{\pair{w}{y}} \in X \times (X \times Y)}
            [\phi(w) = x \mathand \verticalsection{R}{\pair{w}{y}} \in \calI_x]} \\
          & = \preimage{\psi}{\set{\pair{x}{\pair{w}{y}} \in X \times (X \times Y)}
            [\verticalsection{S}{\pair{x}{\pair{w}{y}}} \in \calI_x]},
    \end{align*}
    which is \Borel.
  \end{claimproof}

  If $x \in X$, then $\calI_x \defeq \intersection[n \in \N][\calJ_{\phi_n(x)}]$ is
  a (weakly-ccc-on-\Borel) $\sigma$-ideal for which $\equivalenceclass{x}
  {E} \notin \calI_x$. Clearly $x \mapsto \calI_x$ is $E$-invariant, and Claim
  \ref{idealistic:BorelonBorel} ensures that it is strongly \Borel-on-\Borel.
\end{propositionproof}

\section{Cores} \label{secCores}

Given $n \in \N$ and a set $X$, we use $\sets{n}{X}$ to denote the set of
all subsets of $X$ of cardinality $n$. A \definedterm{partial transversal} of
an equivalence relation $F$ on $X$ is a set $Y \subseteq X$ whose
intersection with each $F$-class consists of at most one point. Given a
superequivalence relation $E$ of $F$, we use $\sets{n}{X}[E][F]$ to denote
the set of all partial transversals $a \in \sets{n}{X}$ of $F$ that are
contained in a single $E$-class. Define $\sets{\le n}{X} \defeq \union[i \le n]
[\sets{i}{X}]$ and $\sets{\le n}{X}[E][F] \defeq \union[i \le n][{\sets{i}{X}[E]
[F]}]$.

If $X$ is a standard \Borel space, then $\sets{n}{X}$ can be viewed as the
quotient of the set of injective sequences in $\functions{n}{X}$ by the
equivalence relation of enumerating the same set, and equipped with the
quotient \Borel structure. The latter is standard: If $\prec$ is a \Borel strict
linear ordering of $X$, then it is easy to see that the quotient mapping
induces a \Borel isomorphism between $\set{x \in \functions{n}{X}}[x(0)
\prec \cdots \prec x(n-1)]$ and $\sets{n}{X}$. We equip $\sets{\le n}{X}$
with the disjoint union \Borel structure. Complexities of subsets of $\sets{n}
{X}[E][F]$ and $\sets{\le n}{X}[E][F]$ will always be considered relative to
the ambient spaces $\sets{n}{X}$ and $\sets{\le n}{X}$, respectively.

We say that $a, b \in \sets{\le n}{X}$ are \definedterm{$F$-disjoint} if
$\equivalenceclass{a}{F} \intersection \equivalenceclass{b}{F} =
\emptyset$. We abuse notation by using $\sets{\le n}{F}$ to denote the
equivalence relation on $\sets{\le n}{X}$ given by $a \mathrel{\sets{\le n}
{F}} b \iff \equivalenceclass{a}{F} = \equivalenceclass{b}{F}$. If $F$ is
\Borel, then so too is $\sets{\le n}{F}$. Finally, we say that a set $C
\subseteq X$ is a \definedterm{core} for a family $\calA \subseteq \sets{\le
n}{X}$ if it intersects every element of $\calA$. In this section, we describe
the circumstances under which suitable subfamilies of $\sets{\le n}{X}$
admit sufficiently small definable cores.

The following result is essentially a special case of \cite[Theorem 12]
{PunctureSets} (although the formalism is quite different and condition (2)
appears in a slightly weaker form in \cite[Theorem 12]{PunctureSets}, the
proof given there yields the result stated below):

\begin{theorem}[\Caicedo--\Clemens--\Conley--\Miller]
  \label{SmallPunctureSets}
  Let $n \ge 1$, $X$ be a standard \Borel space, $F$ be a co-analytic
  equivalence relation on $X$, and $\calA \subseteq \sets{\le n}{X}$ be an
  analytic family of nonempty sets. Then exactly one of the following holds:
  \begin{enumerate}
    \item There is an $F$-invariant core $C \subseteq X$ for $\calA$ on
      which $F$ has only countably many classes.
    \item There is an uncountable \Borel family $\calP \subseteq \calA$
      consisting of pairwise $F$-disjoint sets.
  \end{enumerate}
\end{theorem}

The following result is essentially the analog of \cite[Proposition
2.3.1]{SmoothIdeal} in which bounded finite index is replaced with
countable index, and the idea underlying the proof of the former is
essentially the same as that underlying the proof of the latter:

\begin{proposition} \label{CountableCores}
  Let $n \ge 1$, $X$ be a standard \Borel space, $E$ be an analytic
  equivalence relation on $X$, $F$ be a \Borel equivalence relation on $X$
  contained in $E$, and $\calA \subseteq \sets{\le n}{X}[E][F]$ be an
  analytic family of nonempty sets such that there is a core for $\calA
  \intersection \sets{\le n}{\equivalenceclass{x}{E}}$ that intersects only
  countably many $F$-classes for all $x \in X$. Then there is an
  $F$-invariant \Borel core $C \subseteq X$ for $\calA$ on which $E$ has
  countable index over $F$.
\end{proposition}

\begin{propositionproof}
  We proceed by induction on $n$. For the case $n = 1$, define $A \defeq
  \set{x \in X}[\set{x} \in \calA]$, note that $E \restriction A$ has countable
  index over $F \restriction A$, appeal to Proposition \ref
  {ReflexionCountableIndex} to obtain an $F$-invariant \Borel set $C
  \supseteq A$ on which $E$ has countable index over $F$, and observe
  that $C$ is a core for $\calA$. 

  We now suppose $n \ge 2$. Given $\calF \subseteq \sets{\le n}{X}$ and
  $a \in \sets{\le n}{X}$, we let $\setinterval[F]{a}{\calF} \defeq \set
  {\equivalenceclass{b}{F}}[b \in \calF \mathand a \subseteq
  \equivalenceclass{b}{F}]$. We build, by reverse recursion, analytic families
  $\calA_k \subseteq \sets{\le n}{X}[E][F]$ and $\calA_k' \subseteq \sets{k}
  {X}[E][F]$ for every $k \le n$ and an $F$-invariant \Borel set $B_k
  \subseteq X$ for every $1 \le k \le n$ satisfying the following conditions:
  \begin{enumerate}
    \item For all $k < n$, $\calA_k = \set{a \in \calA_{k+1}}[a \intersection
      B_{k+1} = \emptyset]$.
    \item For all $k \le n$, $\calA_k' =\set{a \in \sets{k}{X}[E][F]}[\cardinality
      {\setinterval[F]{a}{\calA_k}} > \aleph_0]$.
    \item For all $1 \le k \le n$, $B_k$ is a core for $\calA_k'$ on which $E$
      has countable index over $F$.
  \end{enumerate}

  We start with $\calA_n = \calA$. For $k \le n$, condition (2) uniquely
  defines $\calA_k'$ from $\calA_k$. Moreover, since $\sets{k}{X}[E][F]$ is
  analytic and $\calA_k'$ is the set of $a \in \sets{k}{X}[E][F]$ for which
  $\sets{\le n}{F} \restriction \set{b \in \calA_k}[a \subseteq
  \equivalenceclass{b}{F}]$ has uncountably many classes, Corollary \ref
  {ComplexityCountablyManyClasses} ensures that $\calA_k'$ is analytic.
  Similarly, for $k < n$, condition (1) uniquely defines $\calA_k$ from
  $\calA_{k+1}$ and $B_{k+1}$, and ensures that it is analytic. So it only
  remains to describe the construction of the $B_k$'s. From now on, we fix
  $1 \le k \le n$, assume that the $\calA_l$'s and the $\calA_l'$'s have been
  constructed for all $k \le l \le n$ and the $B_l$'s have been constructed
  for $k < l \le n$, and describe the construction of $B_k$. In the case
  when $k = n$, we clearly have $\calA_n' = \emptyset$, so we can take
  $B_n = \emptyset$. Thus we can assume that $1 \le k < n$. We begin with
  several preliminary claims. The first of these, Claim \ref
  {ClaimCoutableExtensions}, is also true when $k = 0$, and will later be
  used in this case.

  \begin{claim} \label{ClaimCoutableExtensions}
    Let $a \in \calA_k'$ and $x \in X$ be such that $a \union \set{x} \in
    \sets{k+1}{X}[E][F]$. Then $\setinterval[F]{a \union \set{x}}{\calA_k}$
    is countable.
  \end{claim}

  \begin{claimproof}
    Suppose $\setinterval[F]{a \union \set{x}}{\calA_k} \neq \emptyset$. Then
    there is $b \in \calA_k$ such that $a \union \set{x} \subseteq
    \equivalenceclass{b}{F}$. By condition (1), we have $b \intersection
    B_{k+1} = \emptyset$. Since $B_{k+1}$ is $F$-invariant, we have
    $\equivalenceclass{b}{F} \intersection B_{k+1} = \emptyset$, hence $(a
    \union \set{x}) \intersection B_{k+1} = \emptyset$. Since $B_{k+1}$ is a
    core for $\calA_{k+1}'$, we deduce that $a \union \set{x} \notin
    \calA_{k+1}'$. Hence, condition (2) ensures that $\setinterval[F]{a \union
    \set{x}}{\calA_{k+1}}$ is countable. Since $\calA_k \subseteq
    \calA_{k+1}$, we deduce that $\setinterval[F]{a \union \set{x}}{\calA_k}$ is
    also countable.
  \end{claimproof}

  \begin{claim} \label{ClaimAvoidingExtensions}
    Let $a \in \calA_k'$ and $M \subseteq \equivalenceclass{a}{E}$ be a set
    intersecting only countably many $F$-classes. Then there exists $b \in
    \calA_k$ such that $a \subseteq \equivalenceclass{b}{F}$ and $M
    \intersection \equivalenceclass{b}{F} \setminus \equivalenceclass{a}{F} =
    \emptyset$.
  \end{claim}

  \begin{claimproof}
    By Claim \ref{ClaimCoutableExtensions}, for all $x \in M \setminus
    \equivalenceclass{a}{F}$, the set $\setinterval[F]{a \union \set{x}}
    {\calA_k}$ is countable. Since this set only depends on
    $\equivalenceclass{x}{F}$ and $M$ intersects only countably many
    $F$-classes, we deduce that $\union[x \in M \setminus
    \equivalenceclass{a}{F}][{\setinterval[F]{a \union \set{x}}{\calA_k}}]$ is
    countable. Since $\setinterval[F]{a}{\calA_k}$ is uncountable, we deduce
    that $\setinterval[F]{a}{\calA_k} \setminus \union[x \in M \setminus
    \equivalenceclass{a}{F}][{\setinterval[F]{a \union \set{x}}{\calA_k}}]$ is
    nonempty. An element of this set has the form $\equivalenceclass{b}{F}$,
    where $b \in \calA_k$, $a \subseteq \equivalenceclass{b}{F}$, and $x
    \notin \equivalenceclass{b}{F}$ for all $x \in M \setminus
    \equivalenceclass{a}{F}$. From this last condition, we deduce that $M
    \intersection \equivalenceclass{b}{F} \setminus \equivalenceclass{a}{F} =
    \emptyset$.
  \end{claimproof}

  \begin{claim} \label{AuxiliaryCore}
    For every $x \in X$, there exists a core for $\calA_k' \intersection \sets{k}
    {\equivalenceclass{x}{E}}$ that intersects only countably many
    $F$-classes.
  \end{claim}

  \begin{claimproof}
    Fix $x \in X$. By the hypotheses of the proposition, there exists a core
    $D \subseteq X$ for $\calA \intersection \sets{\le n}{\equivalenceclass{x}
    {E}}$ that intersects only countably many $F$-classes. Without loss
    of generality, we can assume that $D$ is $F$-invariant. We will show
    that it is also a core for $\calA_k' \intersection \sets{k}{\equivalenceclass
    {x}{E}}$.

    Towards this end, fix $a \in \calA_k' \intersection \sets{k}
    {\equivalenceclass{x}{E}}$. Then $\calF \defeq \set{b \setminus
    \equivalenceclass{a}{F}}[b \in \calA_k \mathand a \subseteq
    \equivalenceclass{b}{F}]$ is an analytic subset of $\sets{\le n}{X}$
    contained in $\sets{\le n}{\equivalenceclass{x}{E}}$. If $\emptyset \in
    \calF$, then there exists $b \in \calA_k$ such that $\equivalenceclass{a}
    {F} = \equivalenceclass{b}{F}$; since $b \in \calA \intersection \sets{\le n}
    {\equivalenceclass{x}{E}}$, we have $b \intersection D \neq \emptyset$,
    and since $D$ is $F$-invariant, we also have $a \intersection D \neq
    \emptyset$. So, from now on, we will assume that $\emptyset \notin
    \calF$. Hence, we can apply Theorem \ref{SmallPunctureSets} to
    $\calF$. There are two cases.

    \case{1}{There is a core $M \subseteq X$ for $\calF$ on which
    $F$ has only countably many classes.} We can assume that $M \subseteq
    \equivalenceclass{x}{E}$. We apply Claim \ref{ClaimAvoidingExtensions}
    to $a$ and $M$, which yields $b \in \calA_k$ such that $a \subseteq
    \equivalenceclass{b}{F}$ and $M \intersection \equivalenceclass{b}{F}
    \setminus \equivalenceclass{a}{F} = \emptyset$. It follows that $b
    \setminus \equivalenceclass{a}{F}$ is in $\calF$ and does not intersect
    $M$, contradicting the fact that $M$ is a core for $\calF$.

    \case{2}{There is an uncountable family $\calP \subseteq \calF$
    consisting of pairwise $F$-disjoint sets.} Since $D$ intersects only
    countably many $F$-classes, we can find $c \in \calP$ such that $c
    \intersection D = \emptyset$. Fix $b \in \calA_k$ such that $a \subseteq
    \saturation{b}{F}$ and $c = b \setminus \saturation{a}{F}$. Then $b \in
    \calA \intersection \sets{\le n}{\equivalenceclass{x}{E}}$; hence $b
    \intersection D \neq \emptyset$. But we have seen that $b \setminus
    \equivalenceclass{a}{F} \intersection D = \emptyset$, so it follows that
    $\equivalenceclass{a}{F} \intersection D \neq \emptyset$. Since $D$ is
    $F$-invariant, we deduce that $a \intersection D \neq \emptyset$.
  \end{claimproof}

  Claim \ref{AuxiliaryCore} allows us to apply the induction hypothesis to
  the family $\calA_k'$; this gives us the desired \Borel set $B_k$ and
  completes the recursive construction. We now let $A \defeq \union
  [\calA_0]$.

  \begin{claim}
    The relation $E$ has countable index over $F$ on $A$.
  \end{claim}

  \begin{claimproof}
    If $\emptyset \notin \calA_0'$, then condition (2) implies that $\setinterval
    [F]{\emptyset}{\calA_0}$ is countable. Observe that $\equivalenceclass
    {A}{F} = \union[{\setinterval[F]{\emptyset}{\calA_0}}]$, so $F$ has only
    countably many classes on $\equivalenceclass{A}{F}$, and we are done.
    We can therefore assume that $\emptyset \in \calA_0'$.

    Fix $x \in X$; we will show that $A \intersection \equivalenceclass{x}{E}$
    intersects only countably many $F$-classes. The hypotheses of the
    proposition give us a core $D \subseteq X$ for $\calA \intersection \sets
    {\le n}{\equivalenceclass{x}{E}}$ that intersects only countably many
    $F$-classes. Since $\emptyset \in \calA_0'$, Claim \ref
    {ClaimCoutableExtensions} ensures that, for every $y \in D$, the set
    $\setinterval[F]{\set{y}}{\calA_0}$ is countable, hence the set $\union
    [{\setinterval[F]{\set{y}}{\calA_0}}]$ intersects only countably many
    $F$-classes. Since this set only depends on $\equivalenceclass{y}{F}$
    and $D$ intersects only countably many $F$-classes, we deduce that
    the set $H \defeq \union[y \in D][{\union[{\setinterval[F]{\set{y}}
    {\calA_0}}]}]$ intersects only countably many $F$-classes.

    It only remains to show that $A \intersection \equivalenceclass{x}{E}
    \subseteq H$. Let $z \in A \intersection
    \equivalenceclass{x}{E}$. Then there exists $a \in \calA_0 \intersection
    \sets{\le n}{\equivalenceclass{x}{E}}$ such that $z \in a$. Since $D$ is a
    core for $\calA_0 \intersection \sets{\le n}{\equivalenceclass{x}{E}}$,
    there exists $y \in a \intersection D$. Then $\equivalenceclass{a}{F}
    \in \setinterval[F]{\set{y}}{\calA_0}$, so $z \in \union[{\setinterval[F]{\set
    {y}}{\calA_0}}]$, thus $z \in H$.
  \end{claimproof}

  We can now complete the proof of the proposition: By Proposition \ref
  {ReflexionCountableIndex}, there exists an $F$-invariant \Borel set $B
  \supseteq A$ on which $E$ has countable index over $F$. This set is a
  core for $\calA_0$. By condition (1), $B_{k+1}$ is a core for $\calA_{k+1}
  \setminus \calA_k$ for every $k < n$. We deduce that $C \defeq B \union
  \left( \union[1 \le k \le n][B_k] \right)$ is a core for $\calA$, and it is clear
  that $C$ is a \Borel set on which $E$ has countable index over $F$.
\end{propositionproof}

Combining Proposition \ref{CountableCores} and Theorem \ref
{SmallPunctureSets}, we obtain:

\begin{corollary} \label{CountableCoresDichotomy}
  Let $n \ge 1$, $X$ be a standard \Borel space, $E$ be an analytic
  equivalence relation on $X$, $F$ be a \Borel equivalence relation on $X$
  contained in $E$, and $\calA \subseteq \sets{\le n}{X}[E][F]$ be an
  analytic family of nonempty sets. Then exactly one of the following holds:
  \begin{enumerate}
    \item There exists an $F$-invariant \Borel core $C \subseteq X$ for
      $\calA$ on which $E$ has countable index over $F$.
    \item There exists an uncountable \Borel set $\calP \subseteq \calA$ of
      pairwise $F$-disjoint subsets of a single $E$-class.
  \end{enumerate}
\end{corollary}

\begin{corollaryproof}
  To see that the two conditions are mutually exclusive, suppose that we
  are in case (2) and fix a core $C$ for $\calA$. Then $C$ contains
  uncountably many pairwise $F$-inequivalent elements of $\union[\calP]$.
  But $\union[\calP]$ is contained in an $E$-class, so $E$ does not have
  countable index over $F$ on $C$.

  Suppose now that we are not in case (2). Then, for every $x \in X$,
  Theorem \ref{SmallPunctureSets} applied to $\calA \intersection \sets
  {\le n}{\equivalenceclass{x}{E}}$ yields a core for $\calA \intersection \sets
  {\le n}{\equivalenceclass{x}{E}}$ intersecting only countably many
  $F$-classes. Hence, we can apply Proposition \ref{CountableCores},
  which ensures that we are in case (1).
\end{corollaryproof}

We now state and prove a technical consequence of Corollary \ref
{CountableCoresDichotomy} that will be useful later. Given sets $X$ and
$Y$, an equivalence relation $F$ on $X \times Y$, and $y \in Y$, we use
$\doublehorizontalsection{F}{y}$ to denote the equivalence relation on $X$
given by $x \mathrel{\doublehorizontalsection{F}{y}} x' \iff \pair{x}{y}
\mathrel{F} \pair{x'}{y}$.

\begin{lemma} \label{CoreMeager}
  Let $X$ and $Y$ be \Polish spaces, $A \subseteq Y$ be an analytic set,
  $E$ be an analytic equivalence relation on $X \times Y$, and $F$ be a
  \Borel equivalence relation on $X \times Y$ contained in $E$ such that:
  \begin{enumerate}
    \renewcommand{\theenumi}{\alph{enumi}}
    \item For every $y \in A$, $X \times \set{y}$ is contained in a single
      $E$-class.
    \item There exists $n \ge 1$ such that, for every $y \in A$, the
      equivalence relation $\doublehorizontalsection{F}{y}$ has at most $n$
      non-meager classes.
  \end{enumerate}
  Then at least one of the following conditions holds:
  \begin{enumerate}
    \item There is an $F$-invariant \Borel set $C \subseteq X \times Y$, on
      which $E$ has countable index over $F$, such that $\horizontalsection
      {C}{y}$ is non-meager for every $y \in A$.
    \item There is a continuous mapping $\psi \from \Cantorspace \to A$
      such that $X \times \image{\psi}{\Cantorspace}$ is contained in a single
      $E$-class and $\preimage{((\id{X} \times \psi) \times (\id{X} \times \psi))}
      {F}$ is meager.
  \end{enumerate}
\end{lemma}

\begin{lemmaproof}
  Let $\calR$ be the set of $\pair{y}{a} \in A \times \sets{\le n}{X \times Y}[E]
  [F]$ with the property that $\equivalenceclass{a}{F} = \set{z \in X \times Y}
  [\horizontalsection{\equivalenceclass{z}{F}}{y} \text{ is non-meager}]$.
  The inclusion $\equivalenceclass{a}{F} \subseteq \set{z \in X \times Y}
  [\horizontalsection{\equivalenceclass{z}{F}}{y} \text{ is non-meager}]$
  can be written as
  \begin{equation*}
    \forall z \in a \existnonmeagerlymany x \in X \ \pair{x}{y} \mathrel{F} z,
  \end{equation*}
  so, by \cite[Theorem 16.1]{Kechris}, the set of $\pair{y}{a} \in Y \times
  \sets{\le n}{X \times Y}[E][F]$ satisfying it is \Borel. Similarly, the reverse
  inclusion can be written as
  \begin{equation*}
    \forcomeagerlymany x \in X \ (\existnonmeagerlymany x' \in X \ \pair{x}
    {y} \mathrel{F} \pair{x'}{y} \implies \exists z \in a \ \pair{x}{y} \mathrel
    {F} z),
  \end{equation*}
  so the set of $\pair{y}{a} \in Y \times \sets{\le n}{X \times Y}[E][F]$
  satisfying it is also \Borel. It follows that $\calR$ is an analytic subset of
  $Y \times \sets{\le n}{X\times Y}$, so $\calA \defeq \image{\projection
  [{\sets{\le n}{X \times Y}}]}{\calR}$ is an analytic subset of $\sets{\le n}{X
  \times Y}[E][F]$.

  Suppose that $\emptyset \in \calA$. Then there exists $y_0 \in A$ such
  that $\pair{y_0}{\emptyset} \in \calR$, that is, $\horizontalsection
  {\equivalenceclass{z}{F}}{y_0}$ is meager for all $z \in X \times Y$. Let
  $\psi \from \Cantorspace \to A$ be the constant mapping with value $y_0$
  and $F' \defeq \preimage{((\id{X} \times \psi) \times (\id{X} \times \psi))}
  {F}$. Then, for all $\pair{x}{u} \in X \times \Cantorspace$ and $u' \in
  \Cantorspace$, we have $\horizontalsection{\equivalenceclass{\pair{x}{u}}
  {F'}}{u'} = \horizontalsection{\equivalenceclass{\pair{x}{y_0}}{F}}{y_0}$, so
  $\horizontalsection{\equivalenceclass{\pair{x}{u}}{F'}}{u'}$ is meager; by
  \Kuratowski--\Ulam's theorem (see, for example, \cite[Theorem 8.41]
  {Kechris}), it follows that $\equivalenceclass{\pair{x}{u}}{F'}$ is meager,
  hence $F'$ is meager. It also follows from condition (a) that $X \times
  \image{\psi}{\Cantorspace}$ is contained in a single $E$-class, so $\psi$
  witnesses that we are in case (2).

  So, from now on, we can assume that $\emptyset \notin \calA$, and
  therefore apply Corollary \ref{CountableCoresDichotomy} to the family
  $\calA$. There are two cases.

  \case{1}{There is an $F$-invariant \Borel core $C \subseteq X
  \times Y$ for $\calA$ on which $E$ has countable index over $F$.} We will
  show that $C$ witnesses that we are in case (1). Let $y \in A$; we will
  show that $\horizontalsection{C}{y}$ is non-meager. By condition (b), the
  equivalence relation $\doublehorizontalsection{F}{y}$ has at most $n$
  non-meager classes; we fix a set $\widetilde{a}$ of representatives for
  these classes and let $a \defeq \widetilde{a} \times \set{y}$. Elements
  of $a$ are clearly pairwise $F$-inequivalent and pairwise $E$-equivalent
  by condition (a), so $a \in \sets{\le n}{X}[E][F]$. Moreover,
  $\equivalenceclass{a}{F}$ is the union of all $F$-classes
  $\equivalenceclass{z}{F}$ for which $\horizontalsection{\equivalenceclass
  {z}{F}}{y}$ is non-meager, so $\pair{y}{a} \in \calR$, thus $a \in \calA$.
  Since $C$ is a core for $\calA$, we can find $x \in X$ such that $\pair{x}
  {y} \in C \intersection a$. It follows that $\horizontalsection
  {\equivalenceclass{\pair{x}{y}}{F}}{y}$ is non-meager, and since $C$ is
  $F$-invariant, we have $\equivalenceclass{\pair{x}{y}}{F} \subseteq C$,
  so $\horizontalsection{C}{y}$ is non-meager.

  \case{2}{There exists an uncountable \Borel set $\calP \subseteq
  \calA$ of pairwise $F$-disjoint subsets of a single $E$-class.} Endow
  $\sets{\le n}{X \times Y}$ with any \Polish topology compatible with its
  standard \Borel structure. We can assume that $\calP$ is homeomorphic
  to Cantor space. By \Jankov--\vonNeumann's uniformization theorem
  (see, for example, \cite[Theorem 18.1]{Kechris}), we can find a \Baire
  measurable mapping $f \from \calP \to Y$ such that, for all $a \in \calP$,
  $\pair{f(a)}{a} \in \calR$. The mapping $f$ is continuous on a comeager
  subset of $\calP$, so, by shrinking $\calP$ if necessary, we can assume
  that $f$ is continuous on $\calP$. Since elements of $\calP$ are
  nonempty, pairwise $F$-disjoint, and contained in the same $E$-class, it
  follows that the sets $\set{z \in X \times Y}[\horizontalsection
  {\equivalenceclass{z}{F}}{f(a)} \text{ is non-meager}]$, for $a \in \calP$,
  are nonempty, pairwise disjoint, and contained in the same $E$-class. In
  particular, the mapping $f$ is one-to-one, so $P \defeq \image{f}{\calP}$
  is homeomorphic to \Cantor space. It is clear that $P \subseteq A$;
  moreover, the sets $\set{z \in X \times Y}[\horizontalsection
  {\equivalenceclass{z}{F}}{y} \text{ is non-meager}]$, for $y \in P$, are
  nonempty, pairwise $F$-disjoint, and contained in the same $E$-class.
  Since, for all $y \in P$, the set $\set{z \in X \times Y}[\horizontalsection
  {\equivalenceclass{z}{F}}{y} \text{ is non-meager}]$ intersects $X \times
  \set{y}$ and $X \times \set{y}$ is contained in a single $E$-class, we
  deduce that $X \times P$ is contained in a single $E$-class.

  We now let $\psi \from \Cantorspace \to P$ be a homeomorphism and
  show that $\psi$ witnesses that we are in case (2), or equivalently, that
  $F \restriction (X \times P)$ is meager. Suppose not. Then, by
  \Kuratowski--\Ulam's theorem, there exist $z \in X \times P$ for which
  $\equivalenceclass{z}{F}$ is non-meager in $X \times P$ and distinct
  $y_0, y_1 \in P$ for which $\horizontalsection{\equivalenceclass{z}{F}}
  {y_0}$ and $\horizontalsection{\equivalenceclass{z}{F}}{y_1}$ are
  non-meager, contradicting the fact that the sets of the form $\set{z \in X
  \times Y}[\horizontalsection{\equivalenceclass{z}{F}}{y} \text{ is
  non-meager}]$, for $y \in P$, are pairwise disjoint.
\end{lemmaproof}

\section{Aligned mappings} \label{secAligned}

In this section and the next, we will often work with spaces of the form
$\CantorCantorspace[m][n]$ where $m, n \in \N \union \set{\N}$. If $m \le
m'$, $n \le n'$, $s \in \CantorCantorspace[m][n]$, and $t \in
\CantorCantorspace[m'][n']$, then we write $s \extendedby t$ to indicate
that $s(j)(i) = t(j)(i)$ for all $i < m$ and $j < n$. Let $\extensions{s} \defeq
\set{x \in \CantorCantorspace[n]}[s \extendedby x]$ for all $m, n \in \N$ and
$s \in \CantorCantorspace[m][n]$. The family $\set{\extensions{s}}[{m \in
\N \mathand s \in \CantorCantorspace[m][n]}]$ forms a basis of clopen
subsets of $\CantorCantorspace[n]$. If $n \in \N$, $m, n' \in \N \union \set
{\N}$, $s \in \CantorCantorspace[m][n]$, and $t \in \CantorCantorspace[m]
[n']$, then we denote by $s \concatenation t$ the \definedterm{horizontal
concatenation} of $s$ and $t$, that is, the element of $\CantorCantorspace
[m][n+n']$ given by $(s \horizontalconcatenation t)(j)(i) = s(j)(i)$ for $\pair{i}
{j} \in m \times n$ and $(s \horizontalconcatenation t)(n + j)(i) = t(j)(i)$ for
$\pair{i}{j} \in m \times n'$. When $n' = 1$, we identify
$\CantorCantorspace[m][1]$ with $\Cantorspace[m]$, so that, for $s \in
\CantorCantorspace[m][n]$ and $t \in \Cantorspace[m]$, the horizontal
concatenation can be written as $s \horizontalconcatenation t$. Similarly,
for $m \in \N$, $m', n \in \N \union \set{\N}$, $s \in \CantorCantorspace[m]
[n]$, and $t \in \CantorCantorspace[m'][n]$, we denote by $s
\verticalconcatenation t$ the \definedterm{vertical concatenation} of $s$
and $t$, that is, the element of $\CantorCantorspace[m+m'][n]$ given by
$(s \verticalconcatenation t)(j)(i) = s(j)(i)$ for $\pair{i}{j} \in m \times n$ and
$(s \verticalconcatenation t)(j)(m+i) = t(j)(i)$ for $\pair{i}{j} \in m' \times n$.

For $f \from \CantorCantorspace[m][n] \to \CantorCantorspace[p][q]$ and
$g \from \CantorCantorspace[m'][n'] \to \CantorCantorspace[p'][q']$, where
$m \le m'$, $n \le n'$, $p \le p'$, and $q \le q'$, we write $f \extendedby g$
when, for all $s \in \CantorCantorspace[m][n]$ and $t \in \CantorCantorspace
[m'][n']$, we have $s \extendedby t \implies f(s) \extendedby g(t)$. In the
special case when $m, n, p, q \in \N$, $f \from \CantorCantorspace[m][n]
\to \CantorCantorspace[p][q]$, and $g \from \CantorCantorspace[n] \to
\CantorCantorspace[q]$, we have $f \extendedby g$ if and only if $\image
{g}{\extensions{s}} \subseteq \extensions{f(s)}$ for all $s \in
\CantorCantorspace[m][n]$.

For $m, n \in \N \union \set{\N}$ and $k \in \N$, let $\Fzero{k}[m][n]$ denote
the equivalence relation on $\CantorCantorspace[m][n]$ given by $s
\mathrel{\Fzero{k}[m][n]} t \iff \forall k \le l < n \ s(l) = t(l)$. Note that this is
equality when $k = 0$ and the complete equivalence relation on
$\CantorCantorspace[m][n]$ when $k \ge n$. In the special case when $m
= n = \N$, we use $\Fzero{k}$ to denote the corresponding equivalence
relation. Observe that $\Eone = \union[k \in \N][\Fzero{k}]$.

For $n \in \N$, an \definedterm{$n$-dimensional aligned mapping} is a
continuous reduction $\phi \from \CantorCantorspace[n] \to
\CantorCantorspace[n]$ of $\sequence{\Fzero{k}[\N][n]}[k < n]$ to itself.
We denote by $\alignedmaps{n}$ the set of all such mappings.

It is clear that $\alignedmaps{n}$ contains the identity and is closed under
composition. In the rest of the paper, $\alignedmaps{n}$ will be viewed as
a subset of $\continuousfunctions{\CantorCantorspace[n]}
{\CantorCantorspace[n]}$ and endowed with the subspace topology.

For all $m, n \in \N$, we let $\restrictedalignedmaps{m}{n} \defeq \set{\phi
\in \alignedmaps{n}}[\id{\CantorCantorspace[m][n]} \extendedby \phi]$. In
other words, a mapping $\phi \from \CantorCantorspace[n] \to
\CantorCantorspace[n]$ belongs to $\restrictedalignedmaps{m}{n}$ if and
only if $\phi \in \alignedmaps{n}$ and $\image{\phi}{\extensions{s}}
\subseteq \extensions{s}$ for all $s \in \CantorCantorspace[m][n]$.

\begin{lemma} \label{BasisNbhdId}
  Fix $n \in \N$. Then $\sequence{\restrictedalignedmaps{m}{n}}[m \in \N]$
  is a basis of neighborhoods of the identity in $\alignedmaps{n}$ and the
  $\restrictedalignedmaps{m}{n}$'s are closed under composition.
\end{lemma}

\begin{lemmaproof}
  The fact that $\restrictedalignedmaps{m}{n}$ contains the identity and is
  closed under composition is clear from the definition. Moreover, we have
  $\restrictedalignedmaps{m}{n} = \alignedmaps{n} \intersection \intersection
  [{s \in \CantorCantorspace[m][n]}][\restrictedcontinuousfunctions
  {\extensions{s}}{\extensions{s}}]$, so $\restrictedalignedmaps{m}{n}$ is
  open.

  To see that $\sequence{\restrictedalignedmaps{m}{n}}[m \in \N]$ is a
  basis of neighborhoods of the identity, take any neighborhood $\calU$ of
  the identity in $\alignedmaps{n}$. By Lemma \ref{BasisCompactOpen},
  we can assume that $\calU = \alignedmaps{n} \intersection \intersection
  [i \in I][\restrictedcontinuousfunctions{K_i}{U_i}]$, where $\sequence{K_i}
  [i \in I]$ is a finite partition of $\CantorCantorspace[n]$ into nonempty
  clopen subsets and $\sequence{U_i}[i \in I]$ is a sequence of open
  subsets of $\CantorCantorspace[n]$. Subdividing the $K_i$'s if necessary,
  we can assume that $I = \CantorCantorspace[m][n]$ for some $m \in \N$,
  and $K_s = \extensions{s}$ for all $s \in \CantorCantorspace[m][n]$.
  Since $\id{\CantorCantorspace[m]} \in \calU$, we have $\extensions{s}
  \subseteq U_s$ for all $s \in \CantorCantorspace[m][n]$; it follows that
  $\restrictedalignedmaps{m}{n} \subseteq \calU$.
\end{lemmaproof}

The following lemma is a special case of \cite[Proposition 2.6]
{TreeableEqRel}:

\begin{lemma} \label{OpenAlignedMappingLemma}
  Let $m, n \in \N$ and $\Omega$ be a family of open subsets of
  $\CantorCantorspace[n]$ which is downwards closed under inclusion and
  has the property that $\union[\Omega]$ is dense in $\CantorCantorspace
  [n]$. Then there exist $m' \ge m$ and a reduction $\psi \from
  \CantorCantorspace[m][n] \to \CantorCantorspace[m'][n]$ of $\sequence
  {\Fzero{k}[m][n]}[k < n]$ to $\sequence{\Fzero{k}[m'][n]}[k < n]$, with $\id
  {\CantorCantorspace[m][n]} \extendedby \psi$, such that $\extensions
  {\psi(s)} \in \Omega$ for all $s \in \CantorCantorspace[m][n]$.
\end{lemma}

\begin{corollary} \label{OpenAlignedMapping}
  Let $m, n \in \N$ and $\Omega$ be a family of open subsets of
  $\CantorCantorspace[n]$ which is downwards closed under inclusion and
  has the property that $\union[\Omega]$ is dense in $\CantorCantorspace
  [n]$. Then there exists an open mapping $\phi \in \restrictedalignedmaps
  {m}{n}$ such that $\image{\phi}{\extensions{s}} \in \Omega$ for all $s \in
  \CantorCantorspace[m][n]$.
\end{corollary}

\begin{corollaryproof}
  Fix $m'$ and $\psi$ as given by Lemma \ref
  {OpenAlignedMappingLemma}. Define $\phi \from \CantorCantorspace[n]
  \to \CantorCantorspace[n]$ by $\phi(s \verticalconcatenation x) \defeq \psi
  (s) \verticalconcatenation x$ for all $s \in \CantorCantorspace[m][n]$ and
  $x \in \CantorCantorspace[n]$. Then $\phi$ satisfies the desired
  conditions.
\end{corollaryproof}

The next lemma gives an example of a situation where families $\Omega$
satisfying the hypotheses of Lemma \ref{OpenAlignedMappingLemma} and
Corollary \ref{OpenAlignedMapping} naturally appear. It will be used in
conjunction with the latter in the proof of our first dichotomy.

\begin{lemma} \label{DecidingCategory}
  Let $X$ be a \Polish space, $E$ be a \Baire measurable equivalence
  relation on $X$, and $\Omega$ be the set of all open sets $U \subseteq
  X$ for which $E \restriction U$ is meager or comeager. Then $\union
  [\Omega]$ is dense in $X$.
\end{lemma}

\begin{lemmaproof}
  Suppose not. Then we can find a nonempty open set $U \subseteq X$
  such that no further nonempty subset of $U$ belongs to $\Omega$. In
  particular, $E \restriction U$ is non-meager, so, by \Kuratowski--\Ulam's
  theorem, there exists $x \in U$ such that the class $\equivalenceclass{x}
  {E} \intersection U$ is \Baire measurable and non-meager. So we can
  find a nonempty open subset $V \subseteq U$ such that
  $\equivalenceclass{x}{E}$ is comeager in $V$. Hence $V \in \Omega$, a
  contradiction.
\end{lemmaproof}

The following result is a version of \Mycielski's theorem for aligned
mappings. It is a particular case of \cite[Proposition 2.10]{TreeableEqRel},
although one should note that the statement of the latter is missing the
hypothesis that $\phi$ is a reduction of $\sequence{\Fzero{k}[m][n]}[k <
n]$ to $\sequence{\Fzero{k}[m'][n]}[k < n]$, which holds in the special
case we require.

\begin{proposition} \label{MycielskiAligned}
  Let $m \in \N$, $n \ge 1$, and $R$ be a comeager binary relation on
  $\CantorCantorspace[n]$. Then there is a homomorphism $\phi \in
  \restrictedalignedmaps{m}{n}$ from $\setcomplement{\Fzero{n-1}[n]}$ to
  $R$.
\end{proposition}

\section{Two dichotomies} \label{secDichoto}

In this section, we prove our two technical dichotomies---Theorems \ref
{UnionCohomomorphism} and \ref{UnionReduction}---and use the
second to prove the \Kechris--\Louveau dichotomy. We start by
recalling the following well-known fact:

\begin{proposition}[see {\cite[Proposition 2.2]{TreeableEqRel}}]
  \label{EoneIndex}
  Fix $k \in \N$ and a \Baire measurable set $B \subseteq
  \CantorCantorspace$. If $\Fzero{k+1}$ has countable index over $\Fzero
  {k}$ on $B$, then $B$ is meager.
\end{proposition}

\begin{theorem} \label{UnionCohomomorphism}
  Let $X$ be a \Polish space, $E$ be an analytic equivalence relation on
  $X$, and $\sequence{E_n}[n \in \N]$ be a sequence of \Borel
  subequivalence relations of $E$. Then exactly one of the following holds:
  \begin{enumerate}
    \item There exists a cover $\sequence{B_n}[n \in \N]$ of $X$ with the
      property that $B_n$ is an $E_n$-invariant \Borel set on which $E$ has
      countable index over $E_n$ for all $n \in \N$.
    \item There is a continuous homomorphism $\phi \from
      \CantorCantorspace \to X$ from $\sequence{\Eone \setminus \Fzero
      {n}}[n \in \N]$ to $\sequence{E \setminus \union[i \le n][E_i]}[n \in \N]$.
\end{enumerate}
\end{theorem}

\begin{theoremproof}
  We first show that the two conditions are mutually exclusive. Suppose,
  towards a contradiction, that both hold. For every $n \in \N$, define $B_n'
  \defeq \preimage{\phi}{B_n}$. Then $\sequence{B_n'}[n \in \N]$ is a
  covering of $\CantorCantorspace$, so there exists $n_0 \in \N$ such that
  $B_{n_0}'$ is non-meager. By Proposition \ref{EoneIndex}, $\Fzero
  {n_0+1}$ does not have countable index over $\Fzero{n_0}$ on
  $B_{n_0}'$, so neither does $\Eone$. As $\phi \restriction B_{n_0}'$ is
  a homomorphism from $(\Eone \setminus \Fzero{n_0}) \restriction
  B_{n_0}'$ to $(E \setminus E_{n_0}) \restriction B_{n_0}$, it follows that
  $E$ does not have countable index over $E_{n_0}$ on $B_{n_0}$, the
  desired contradiction.

  We now show that at least one of the two conditions holds. We begin by
  fixing some notation, definitions, and conventions. For every $k \in \N$,
  let $R_k \defeq \union[i \le k][E_i]$ and fix a \Polish space $Y_k$ and a
  continuous surjection $\pi_k \from Y_k \to E \setminus R_k$.

  Note that the mapping $\pair{\pair{x}{x'}}{y} \mapsto \pair{x
  \horizontalconcatenation y}{x' \horizontalconcatenation y}$ is a
  homeomorphism from $\setcomplement{\Fzero{k}[k+1]} \times
  \CantorCantorspace[n-k-1]$ to $\Fzero{k+1}[n] \setminus \Fzero{k}[n]$ for
  all $n \in \N \union \set{\N}$ and $k < n$. In the rest of the proof, we will
  identify these spaces via this homeomorphism. As a consequence, for
  instance, if $k < m \le n$ and $g \from \Fzero{k+1}[n] \setminus \Fzero{k}
  [n] \to Y_k$, then it is consistent with our earlier notation to use
  $\horizontalsection{g}{y}$ to denote the mapping $\horizontalsection{g}{y}
  \from \Fzero{k+1}[m] \setminus \Fzero{k}[m] \to Y_k$ given by
  $\horizontalsection{g}{y}(x, x') \defeq g(x \horizontalconcatenation y, x'
  \horizontalconcatenation y)$ for all $\pair{x}{x'} \in \Fzero{k+1}[m]
  \setminus \Fzero{k}[m]$ and $y \in \CantorCantorspace[n - m]$.

  Suppose that $n \in \N \union \set{\N}$ and $m \le n$.
  \begin{itemize}
    \item For $f \in \continuousfunctions{\CantorCantorspace[m]}{X}$,
      define $\sectionpower{f}{n} \in \continuousfunctions
      {\CantorCantorspace[n]}{X}$ by $\sectionpower{f}{n}(x) \defeq f(x
      \restriction m)$.
    \item For $\calA \subseteq \continuousfunctions{\CantorCantorspace[m]}
      {X}$, we abuse notation by using $\sectionpower{\calA}{n}$ to denote
      the set of all $f \in \continuousfunctions{\CantorCantorspace[n]}{X}$
      such that $\horizontalsection{f}{y} \in \calA$ for all $y \in
      \CantorCantorspace[n-m]$. Note that if $f \in \calA$, then
      $\sectionpower{f}{n} \in \sectionpower{\calA}{n}$.
    \item For $k < m$ and $g \in \continuousfunctions{\Fzero{k+1}[m]
      \setminus \Fzero{k}[m]}{Y_k}$, define $\sectionpower{g}{n} \in
      \continuousfunctions{\Fzero{k+1}[n] \setminus \Fzero{k}[n]}{Y_k}$ by
      $\sectionpower{g}{n}(x, x') \defeq g(x \restriction m, x' \restriction m)$.
    \item For $k < m$ and $\calA \subseteq \continuousfunctions{\Fzero
      {k+1}[m] \setminus \Fzero{k}[m]}{Y_k}$, we abuse notation by using
      $\sectionpower{\calA}{n}$ to denote the set of $g \in
      \continuousfunctions{\Fzero{k+1}[n] \setminus \Fzero{k}[n]}{Y_k}$ such
      that $\horizontalsection{g}{y} \in \calA$ for all $y \in \CantorCantorspace
      [n-m]$. Note that if $g \in \calA$, then $\sectionpower{g}{n} \in
      \sectionpower{\calA}{n}$.
  \end{itemize}

  By Lemma \ref{CompatibleMetrics}, we can find compatible complete
  metrics $d^n$ on $\continuousfunctions{\CantorCantorspace[n]}{X}$ such
  that $\diam^\N(\sectionpower{\calA}{\N}) \le \diam^n(\calA)$ for all $n \in
  \N$ and $\calA \subseteq \continuousfunctions{\CantorCantorspace[n]}
  {X}$. Similarly, for all $k \in \N$, we can find compatible complete metrics
  $d_k^n$ on $\continuousfunctions{\Fzero{k+1}[n] \setminus \Fzero{k}[n]}
  {Y_k}$ such that $\diam_k^\N(\sectionpower{\calA}{\N}) \le \diam_k^n
  (\calA)$ for all $n > k$ and $\calA \subseteq \continuousfunctions{\Fzero
  {k+1}[n] \setminus \Fzero{k}[n]}{Y_k}$. As it should not lead to confusion,
  we will use $d$ to denote these metrics and $\diam$ to denote the
  corresponding diameters.

  Using Corollary \ref{BasisRightStable} and the fact that $\Fzero{k+1}[n]
  \setminus \Fzero{k}[n]$ can be identified with an open subset of the
  compact zero-dimensional \Polish space $(\CantorCantorspace[k+1]
  \times \CantorCantorspace[k+1]) \times \CantorCantorspace[n-k-1]$, we
  can fix countable bases of nonempty open subsets of
  $\continuousfunctions{\CantorCantorspace[n]}{X}$ and
  $\continuousfunctions{\Fzero{k+1}[n] \setminus \Fzero{k}[n]}{Y_k}$ whose
  elements are right stable for all $n \in \N$ and $k < n$. These bases will
  not be given a name, but we will refer to them by talking about
  \definedterm{basic open subsets} of these spaces.

  An \definedterm{approximation} is a triple of the form $a \defeq \triple{n^a}
  {\calU^a}{\sequence{\calV_k^a}[k < n^a]}$, where $n^a \in \N$, $\calU^a
  \subseteq \continuousfunctions{\CantorCantorspace[n^a]}{X}$ is a basic
  open set with $\diam(\calU^a) \le 1/n^a$, and $\calV_k^a \subseteq
  \continuousfunctions{\Fzero{k+1}[n^a] \setminus \Fzero{k}[n^a]}{Y_k}$ is
  a basic open set with $\diam(\calV_k^a) \le 1 / n^a$ for all $k < n^a$.
  Given approximations $a$ and $b$, we say that $b$ \definedterm
  {extends} $a$ if $n^a < n^b$, $\calU^b \subseteq \sectionpower{(\calU^a)}
  {n^b}$, and $\calV_k^b \subseteq \sectionpower{(\calV_k^a)}{n^b}$ for all
  $k < n^a$. We say that $b$ is an \definedterm{immediate successor} of
  $a$ if it extends $a$ and $n^b = n^a + 1$.

  A \definedterm{configuration} is a triple of the form $\gamma \defeq \triple
  {n^\gamma}{f^\gamma}{\sequence{g_k^\gamma}[k < n^\gamma]}$, where
  $n^\gamma \in \N$, $f^\gamma \from \CantorCantorspace[n^\gamma] \to
  X$ is continuous, and $g_k^\gamma \from \Fzero{k+1}[n^\gamma]
  \setminus \Fzero{k}[n^\gamma]$ $\to Y_k$ is a continuous mapping such
  that $\pair{f^\gamma(x)}{f^\gamma(x')} = \pi_k(g_k^\gamma(x, x'))$ for all
  $\pair{x}{x'} \in \Fzero{k+1}[n^\gamma] \setminus \Fzero{k}[n^\gamma]$
  and $k < n^\gamma$. An immediate consequence of this definition is that
  if $\gamma$ is a configuration, then $f^\gamma$ is a homomorphism from
  $\sequence{\Fzero{k+1}[n^\gamma] \setminus \Fzero{k}[n^\gamma]}[k <
  n^\gamma]$ to $\sequence{E \setminus R_k}[k < n^\gamma]$; in
  particular, it is injective and its range is contained in a single $E$-class.
  We say that a configuration $\gamma$ is \definedterm{compatible} with
  an approximation $a$ if $n^\gamma = n^a$, $f^\gamma \in \calU^a$, and
  $g_k^\gamma \in \calV_k^a$ for all $k < n^\gamma$. We say that a
  configuration $\gamma$ is \definedterm{generically compatible} with a
  set $B \subseteq X$ if $\preimage{(f^\gamma)}{B}$ is comeager in
  $\CantorCantorspace[n^\gamma]$.

  For all $n \in \N$, the set of all configurations $\gamma$ such that
  $n^\gamma = n$ can be identified with a subset $\configurations{n}$ of
  the space $\continuousfunctions{\CantorCantorspace[n]}{X} \times
  \product[k < n][\continuousfunctions{\Fzero{k+1}[n] \setminus \Fzero{k}[n]}
  {Y_k}]$. The latter space is \Polish by Proposition \ref
  {PolishCompactOpen}, and Proposition \ref{ContEval} ensures that
  $\configurations{n}$ is closed, thus \Polish. For an approximation $a$
  and $B \subseteq X$, define $\compatible{a}{B}$ as the set of all
  configurations that are compatible with $a$ and generically compatible
  with $B$. By definition, being compatible with an approximation is an open
  condition. Moreover, if $B$ is \Borel, then it follows from \cite[Theorem
  16.1]{Kechris} and Proposition \ref{ContEval} that the set of all $\gamma
  \in \configurations{n^a}$ that are generically compatible with $B$ is \Borel.
  Hence, $\compatible{a}{B}$ is a \Borel subset of $\configurations{n^a}$.

  For a configuration $\gamma$ and $\phi \in \alignedmaps{n^\gamma}$,
  we abuse notation by using $\gamma \composition \phi$ to denote the
  triple $\delta \defeq \triple{n^\delta}{f^\delta}{\sequence{g_k^\delta}[k <
  n^\delta]}$, where $n^\delta \defeq n^\gamma$, $f^\delta \defeq f^\gamma
  \composition \phi$, and $g_k^\delta \defeq g_k^\gamma \composition
  (\phi \times \phi) \restriction \Fzero{k+1}[n^\gamma] \setminus \Fzero{k}
  [n^\gamma]$ for all $k < n^\delta$. It follows from the definition of an
  aligned mapping that $\gamma \composition \phi$ is a well-defined
  configuration. Given configurations $\gamma, \delta$ and $m \in \N$, we
  write $\delta \preccurlyeq_m \gamma$ if there exists $\phi \in
  \restrictedalignedmaps{m}{n^\gamma}$ such that $\delta = \gamma
  \composition \phi$. It follows from Lemma \ref{BasisNbhdId} that
  $\preccurlyeq_m$ is a quasi-ordering of the set of all configurations.

  \begin{claim} \label{ClaimDownwardsClosed}
    Let $a$ be an approximation. Then there exists $m \in \N$ such that the
    set $\compatible{a}{X}$ of all configurations that are compatible with
    $a$ is downwards closed under $\preccurlyeq_m$.
  \end{claim}

  \begin{claimproof}
    Since $\calU^a$ is a basic open set of $\continuousfunctions
    {\CantorCantorspace[n^a]}{X}$, hence right stable, we can find $p \in \N$
    with the property that $\calU^a \composition \restrictedalignedmaps{p}
    {n^a} = \calU^a$. Similarly, for every $k < n^a$, since $\calV_k^a$ is a
    right-stable open subset of $\continuousfunctions{\Fzero{k+1}[n^a]
    \setminus \Fzero{k}[n^a]}{Y_k}$, there is a neighborhood $\calW_k$ of
    the identity in $\continuousfunctions{\Fzero{k+1}[n^a] \setminus \Fzero
    {k}[n^a]}{\Fzero{k+1}[n^a] \setminus \Fzero{k}[n^a]}$ such that
    $\calV_k^a \composition \calW_k = \calV_k^a$. By Propositions \ref
    {ProductCompactOpen}, \ref{ContCompo}, and \ref
    {SubspaceRangeCompactOpen}, the mapping $\alignedmaps{n} \to
    \continuousfunctions{\Fzero{k+1}[n^a] \setminus \Fzero{k}[n^a]}{\Fzero
    {k+1}[n^a] \setminus \Fzero{k}[n^a]}$, given by $\phi \mapsto (\phi
    \times \phi) \restriction \Fzero{k+1}[n^a] \setminus \Fzero{k}[n^a]$, is
    continuous; hence we can find $m_k \in \N$ such that $(\phi \times \phi)
    \restriction \Fzero{k+1}[n^a] \setminus \Fzero{k}[n^a]$ is in $\calW_k$ for
    all $\phi \in \restrictedalignedmaps{m_k}{n^a}$. Then $m \defeq \max(p,
    m_0, \ldots, m_{n^a - 1})$ is as desired.
  \end{claimproof}

  \begin{claim} \label{MeagerPullback}
    Let $n \ge 1$ and $\Gamma \from \Cantorspace \to \configurations{n}$
    be a continuous mapping and define $F \from \CantorCantorspace[n+1]
    \to X$ by $F(x \horizontalconcatenation y) \defeq f^{\Gamma(y)}(x)$ for
    all $x \in \CantorCantorspace[n]$ and $y \in \Cantorspace$. Then
    $\preimage{(F \times F)}{R_{n-1}}$ is meager in $\CantorCantorspace
    [n+1] \times \CantorCantorspace[n+1]$.
  \end{claim}

  \begin{claimproof}
    Suppose not. Then there exists $i < n$ such that $\preimage{(F \times
    F)}{E_i}$ is non-meager. By Proposition \ref{ContEval}, the
    mapping $F$ is continuous, so $\preimage{(F \times F)}{E_i}$ is \Borel.
    Hence, by \Kuratowski--\Ulam's theorem, we can find $\triple{u}{x}{z} \in
    \CantorCantorspace[n+1] \times \CantorCantorspace[n-1] \times
    \Cantorspace$ such that the set $A \defeq \set{y \in \Cantorspace}[F(u)
    \mathrel{E_i} F(x \horizontalconcatenation y \horizontalconcatenation
    z)]$ is non-meager in $\Cantorspace$. Let $y, y' \in A$ be distinct. Then
    $\neg x \horizontalconcatenation y \mathrel{\Fzero{n-1}} x
    \horizontalconcatenation y'$, so $f^{\Gamma(z)}(x
    \horizontalconcatenation y) \mathrel{E \setminus R_{n-1}} f^{\Gamma
    (z)}(x \horizontalconcatenation y)$, thus $\neg F(x
    \horizontalconcatenation y \horizontalconcatenation z) \mathrel{E_i} F(x
    \horizontalconcatenation y' \horizontalconcatenation z)$. Since $E_i$ is
    an equivalence relation, this contradicts the fact that both $F(x
    \horizontalconcatenation y \horizontalconcatenation z)$ and $F(x
    \horizontalconcatenation y' \horizontalconcatenation z)$ are
    $E_i$-related to $F(u)$.
  \end{claimproof}

  We will recursively build a decreasing sequence $\sequence{X_\alpha}
  [\alpha < \omega_1]$ of \Borel subsets of $X$. We start with $X_0 \defeq
  X$, and for limit ordinals $\lambda$, we let $X_\lambda \defeq
  \intersection[\alpha < \lambda][X_\alpha]$. We now fix $\alpha <
  \omega_1$ and assume $X_\alpha$ has been constructed; we describe
  how to construct $X_{\alpha + 1}$. We denote by $S_\alpha$ the set of
  all approximations $a$ for which $\compatible{a}{X_\alpha} \neq
  \emptyset$. We say that an approximation $a$ is \definedterm
  {$\alpha$-terminal} if $a$ has no immediate successor in $S_\alpha$; we
  denote by $T_\alpha$ the set of all $\alpha$-terminal approximations.

  \begin{claim} \label{CaracTerminal}
    Let $a$ be an $\alpha$-terminal approximation. Then there exists an
    $E_{n^a}$-invariant \Borel set $B \subseteq X$, on which $E$ has
    countable index over $E_{n^a}$, such that $\compatible{a}{X_\alpha
    \setminus B} = \emptyset$.
  \end{claim}

  \begin{claimproof}
    Proposition \ref{ContEval} ensures that the function $\evaluation \from
    \CantorCantorspace[n^a] \times \configurations{n^a} \to X$, given by
    $\evaluation(x,\gamma) \defeq f^\gamma(x)$, is continuous. Fix $m \in
    \N$ as given by Claim \ref{ClaimDownwardsClosed} applied to the
    approximation $a$. Let $\calB$ be the set of all $\gamma \in \compatible
    {a}{X_\alpha}$ such that $\preimage{(f^\gamma \times f^\gamma)}
    {E_{n^a}} \restriction \extensions{s}$ is either meager or comeager for all
    $s \in \CantorCantorspace[m][n^a]$. The continuity of $\evaluation$ and
    \cite[Theorem 16.1]{Kechris} ensures that $\calB$ is a \Borel subset of
    $\configurations{n^a}$. For all $\gamma \in \calB$ and $s \in
    \CantorCantorspace[m][n^a]$, \Kuratowski--\Ulam's theorem implies that
    each class of $\preimage{(f^\gamma \times f^\gamma)}{E_{n^a}}$ is
    either meager or comeager in $\extensions{s}$, so $\preimage
    {(f^\gamma \times f^\gamma)}{E_{n^a}}$ has at most $2^{m n^a}$
    non-meager classes. Along with the fact that $\image{f^\gamma}
    {\CantorCantorspace[n^a]}$ is contained in a single $E$-class for all
    $\gamma \in \calB$, this ensures that we can apply Lemma \ref
    {CoreMeager} to the spaces $X' \defeq \CantorCantorspace[n^a]$ and
    $Y' \defeq \configurations{n^a}$, the set $\calB \subseteq Y'$, and the
    equivalence relations $E' \defeq \preimage{(\evaluation \times
    \evaluation)}{E}$ and $E_{n^a}' \defeq \preimage{(\evaluation \times
    \evaluation)}{E_{n^a}}$. There are two cases.

    \case{1}{There is an $E_{n^a}'$-invariant \Borel set $B' \subseteq
    \CantorCantorspace[n^a] \times \configurations{n^a}$, on which $E'$
    has countable index over $E_{n^a}'$, such that $(B')^\gamma$ is
    non-meager in $\CantorCantorspace[n^a]$ for all $\gamma \in
    \calB$.} Lemma \ref{PushForwardCountableIndex} ensures that $E$
    has countable index over $E_{n^a}$ on $\image{\evaluation}{B'}$,
    and since the latter set is analytic, Proposition \ref
    {ReflexionCountableIndex} yields an $E_{n^a}$-invariant \Borel set $B
    \supseteq \image{\evaluation}{B'}$ on which $E$ has countable index
    over $E_{n^a}$. Observe that $(B')^\gamma \subseteq \preimage
    {(f^\gamma)}{B}$ for all $\gamma \in \calB$, so $\preimage{(f^\gamma)}
    {B}$ is non-meager, thus no $\gamma \in \calB$ is generically
    compatible with $X_\alpha \setminus B$. It remains to show that this
    holds of every configuration $\gamma$ that is compatible with $a$.

    We can assume, without loss of generality, that $\gamma \in \compatible
    {a}{X_\alpha}$. Let $\Omega$ be the set of all open subsets $U
    \subseteq \CantorCantorspace[n^a]$ with the property that $\preimage
    {(f^\gamma \times f^\gamma)}{E_{n^a}} \restriction U$ is either meager
    or comeager. By Lemma \ref{DecidingCategory}, $\union[\Omega]$ is
    dense in $\CantorCantorspace[n^a]$. Hence, we can apply Corollary
    \ref{OpenAlignedMapping} to find an open mapping $\phi \in
    \restrictedalignedmaps{m}{n^a}$ such that $\image{\phi}{\extensions{s}}
    \in \Omega$ for all $s \in \CantorCantorspace[m][n^a]$. Now let $\delta
    \defeq \gamma \composition \phi$. Then $\delta \preccurlyeq_m
    \gamma$; hence, by the choice of $m$, $\delta$ is compatible with $a$.
    Moreover, since $\gamma$ is generically compatible with $X_\alpha$
    and $\phi$ is open and one-to-one, it follows that $\delta$ is
    generically compatible with $X_\alpha$, hence $\delta \in \compatible
    {a}{X_\alpha}$. The conditions on $\phi$ imply that if $s \in
    \CantorCantorspace[m][n^a]$, then $\preimage{(f^\gamma \times
    f^\gamma)}{E_{n^a}} \restriction \image{\phi}{\extensions{s}}$ is either
    meager or comeager, hence $\preimage{(f^\delta \times f^\delta)}
    {E_{n^a}} \restriction \extensions{s}$ is either meager or comeager. It
    follows that $\delta \in \calB$. Hence, $\preimage{(f^\delta)}{B}$ is
    non-meager, and since $\phi$ is one-to-one and open, we deduce that
    $\preimage{(f^\gamma)}{B}$ is non-meager, thus $\gamma$ is not
    generically compatible with $X_\alpha \setminus B$.

    \case{2}{There exists a continuous mapping $\Gamma \from
    \Cantorspace \to \calB$ such that $\CantorCantorspace[n^a] \times
    \image{\Gamma}{\Cantorspace}$ is contained in a single $E'$-class and
    the equivalence relation $\preimage{((\id{\CantorCantorspace[n^a]}
    \times \Gamma) \times (\id{\CantorCantorspace[n^a]} \times \Gamma))}
    {E_{n^a}'}$ is meager.} We will show that the approximation $a$ is not
    $\alpha$-terminal. Define $F \defeq \evaluation \composition (\id
    {\CantorCantorspace[n^a]} \times \Gamma) \from \CantorCantorspace
    [n^a+1] \to X$, so that $F(x \horizontalconcatenation y) = f^{\Gamma
    (y)}(x)$ for all $x \in \CantorCantorspace[n^a]$ and $y \in \Cantorspace$.
    Our assumptions on $\Gamma$ imply that $\image{F}
    {\CantorCantorspace[n^a+1]}$ is contained in a single $E$-class and
    $\preimage{(F \times F)}{E_{n^a}}$ is meager. If $n^a \ge 1$,  then
    Claim \ref{MeagerPullback} implies that $\preimage{(F \times F)}
    {R_{n^a - 1}}$ is meager, in which case $\preimage{(F \times F)}
    {R_{n^a}}$ is meager, and this obviously remains true when $n^a = 0$.
    For all $y \in \Cantorspace$, $\Gamma(y)$ is generically compatible with
    $X_\alpha$, so $\forall y \in \Cantorspace \forcomeagerlymany x \in
    \CantorCantorspace[n^a] \ F(x \horizontalconcatenation y) \in X_\alpha$,
    thus \Kuratowski--\Ulam's theorem ensures that $\preimage{F}
    {X_\alpha}$ is comeager in $\CantorCantorspace[n^a+1]$.

    Let $C$ be a dense $G_\delta$ subset of $\CantorCantorspace[n^a+1]
    \times \CantorCantorspace[n^a+1]$ contained in $(\preimage{F}
    {X_\alpha} \times \preimage{F}{X_\alpha}) \setminus \preimage{(F \times
    F)}{R_{n^a}}$. For all $\pair{x}{x'} \in C$, we have $\pair{F(x)}{F(x')} \in
    E \setminus R_{n^a}$, so there exists $u \in Y_{n^a}$ such that
    $\pi_{n^a}(u) = \pair{F(x)}{F(x')}$. Hence, by \Jankov--\vonNeumann's
    uniformisation theorem, we can find a \Baire measurable mapping $G
    \from C \to Y_{n^a}$ such that $\pi_{n^a}(G(x,x')) = \pair{F(x)}{F(x')}$ for
    all $\pair{x}{x'} \in C$. Since \Baire measurable mappings are continuous
    on a comeager set, by shrinking $C$ if necessary, we can assume that
    $G$ is continuous.

    By Proposition \ref{MycielskiAligned}, there exists a homomorphism
    $\phi \in \restrictedalignedmaps{m}{n^a+1}$ from $\setcomplement
    {\Fzero{n^a}[n^a+1]}$ to $C$. Since $\phi$ is aligned, it follows that if $y
    \in \Cantorspace$, then $\phi(x \horizontalconcatenation y)(n^a)$ does
    not depend upon $x$. We denote this value by $\psi(y)$, thereby
    obtaining a continuous mapping $\psi \from \Cantorspace \to
    \Cantorspace$. For all $x \in \CantorCantorspace[n^a]$ and $y \in
    \Cantorspace$, we hence write $\phi(x \horizontalconcatenation y) =
    (\phi^y \restriction n^a)(x) \horizontalconcatenation \psi(y)$, where we
    abuse notation by using $\phi^y \restriction n^a$ to denote the aligned
    mapping $\CantorCantorspace[n^a] \to \CantorCantorspace[n^a]$ given
    by $(\phi^y \restriction n^a)(x) \defeq \phi(x \horizontalconcatenation y)
    \restriction n^a$. To see that $\phi^y \restriction n^a \in
    \restrictedalignedmaps{m}{n^a}$, note that if $s \in \CantorCantorspace
    [m][n^a]$ and $x \in \extensions{s}$, then $x \horizontalconcatenation y
    \in \extensions{s \horizontalconcatenation y \restriction m}$, so $\phi(x
    \horizontalconcatenation y) \in \extensions{s \horizontalconcatenation y
    \restriction m}$, thus $(\phi^y \restriction n^a)(x) \in \extensions{s}$.

    Define $\gamma \defeq \triple{n^\gamma}{f^\gamma}{\sequence
    {g_k^\gamma}[k < n^\gamma]}$, where $n^\gamma \defeq n^a + 1$,
    $f^\gamma \from \CantorCantorspace[n^\gamma] \to X$ is given by
    $f^\gamma \defeq F \composition \phi$, $g_k^\gamma \from \Fzero{k+1}
    [n^\gamma] \setminus \Fzero{k}[n^\gamma] \to Y_k$ is given by
    $g_k^\gamma(x \horizontalconcatenation y, x' \horizontalconcatenation
    y) \defeq g_k^{\Gamma(\psi(y))}((\phi^y \restriction n^a)(x), (\phi^y
    \restriction n^a)(x'))$ for all $k < n^a$, $\pair{x}{x'} \in \Fzero{k+1}[n^a]
    \setminus \Fzero{k}[n^a]$, and $y \in \Cantorspace$, and $g_{n^a}^
    \gamma \defeq G \composition (\phi \times \phi) \restriction
    \setcomplement{\Fzero{n^a}[n^a+1]}$. The continuity of $F$ yields that
    of $f^\gamma$, Proposition \ref{ContEval} ensures that $g_k^\gamma$
    is continuous for all $k < n^a$, and the continuity of $G$ on $\image
    {(\phi \times \phi)}{\setcomplement{\Fzero{n^a}[n^a+1]}}$ yields the
    continuity of $g_{n^a}^\gamma$. Moreover, if $k < n^a$, $\pair{x}{x'} \in
    \Fzero{k+1}[n^a] \setminus \Fzero{k}[n^a]$, and $y \in \Cantorspace$,
    then
    \begin{align*}
      \pi_k
        & (g_k^\gamma(x \horizontalconcatenation y, x'
          \horizontalconcatenation y)) \\
        & = \pi_k(g_k^{\Gamma(\psi(y))}((\phi^y \restriction n^a)(x), (\phi^y
          \restriction n^a)(x'))) \\
        & = (f^{\Gamma(\psi(y))}((\phi^y \restriction n^a)(x)),f^{\Gamma(\psi(y))}
          ((\phi^y \restriction n^a)(x'))) \\
        & = (F((\phi^y \restriction n^a)(x) \horizontalconcatenation \psi(y)),
          F((\phi^y \restriction n^a)(x') \horizontalconcatenation \psi(y))) \\
        & = (F(\phi(x \horizontalconcatenation y)), F(\phi(x'
          \horizontalconcatenation y))) \\
        & = (f^\gamma(x \horizontalconcatenation y), f^\gamma(x'
          \horizontalconcatenation y)).
    \end{align*}
    And if $\pair{x}{x'} \in \setcomplement{\Fzero{n^a}[n^\gamma]}$, then the
    definition of $G$ ensures that $\pi_{n^a}(G(\phi(x), \phi(x'))) = \pair{F(\phi
    (x))}{F(\phi(x'))}$, that is, $\pi_{n^a}(g_{n^a}^\gamma(x, x')) = \pair
    {f^\gamma(x)}{f^\gamma(x')}$. Hence, $\gamma$ is a configuration.
    Moreover, since $\phi$ is a homomorphism from $\setcomplement
    {\Fzero{n^a}[n^a+1]}$ to $C$ and $C \subseteq \preimage{F}{X_\alpha}
    \times \preimage{F}{X_\alpha}$, it follows that $f^\gamma$ takes values
    in $X_\alpha$; in particular, $\gamma$ is generically compatible with
    $X_\alpha$.

    For all $y \in \Cantorspace$, define $\gamma^y \defeq \triple{n^a}
    {(f^\gamma)^y}{\sequence{(g_k^\gamma)^y}[k < n^a]}$. The formulas
    defining $\gamma$ ensure that $\gamma^y = \Gamma(\psi(y))
    \composition (\phi^y \restriction n^a)$; since $\phi^y \restriction n^a \in
    \restrictedalignedmaps{m}{n^a}$, it follows that $\gamma^y$ is a
    configuration and $\gamma^y \preccurlyeq_m \Gamma(\psi(y))$. Since
    $\Gamma(\psi(y))$ is compatible with $a$, our choice of $m$ ensures
    that $\gamma^y$ is compatible with $a$. In particular, if $y \in
    \Cantorspace$, then $(f^{\gamma})^y \in \calU^a$, so $f^\gamma
    \in \sectionpower{(\calU^a)}{n^\gamma}$; similarly, if $k < n^a$ and
    $y \in \Cantorspace$, then $(g_k^\gamma)^y \in \calV_k^a$, so
    $g_k^\gamma \in \sectionpower{(\calV_k^a)}{n^\gamma}$. Recall that,
    by Lemma \ref{ArrowPreservesOpen}, $\sectionpower{(\calU^a)}
    {n^\gamma}$ is an open subset of $\continuousfunctions
    {\CantorCantorspace[n^\gamma]}{X}$ and $\sectionpower{(\calV_k^a)}
    {n^\gamma}$ is an open subset of $\continuousfunctions{\Fzero{k+1}
    [n^\gamma] \setminus \Fzero{k}[n^\gamma]}{Y_k}$ for all $k < n^a$.
    Thus, we can find a basic open set $\calU^b \subseteq
    \continuousfunctions{\CantorCantorspace[n^\gamma]}{X}$ with $\diam
    (\calU^b) \le 1 / n^\gamma$ and $f^\gamma \in \calU^b \subseteq
    \sectionpower{(\calU^a)}{n^\gamma}$ and basic open sets $\calV_k^b
    \subseteq \continuousfunctions{\Fzero{k+1}[n^\gamma] \setminus \Fzero
    {k}[n^\gamma]}{Y_k}$ with $\diam(\calV_k^b) \le 1 / n^\gamma$ and
    $g_k^\gamma \in \calV_k^b \subseteq \sectionpower{(\calV_k^a)}
    {n^\gamma}$ for all $k < n^a$. We can also find a basic open set
    $\calV_{n^a}^b \subseteq \continuousfunctions{\setcomplement{\Fzero
    {n^a}[n^\gamma]}}{Y_{n^a}}$ with $\diam(\calV_{n^a}^b) \le 1/n^\gamma$
    and $g_{n^a}^\gamma \in \calV_{n^a}^b$. Letting $n^b \defeq
    n^\gamma$, it follows that $b \defeq \triple{n^b}{\calU^b}{\sequence
    {\calV_k^b}[k < n^b]}$ is an approximation which is an immediate
    successor of $a$ and with which $\gamma$ is compatible. In particular,
    $\compatible{b}{X_\alpha} \neq \emptyset$, so $b \in S_\alpha$, thus
    $a$ is not terminal.
  \end{claimproof}

  For all $a \in T_\alpha$, let $B_\alpha^a$ be a \Borel set as given by
  Claim \ref{CaracTerminal}. Define $X_{\alpha + 1} \defeq X_\alpha
  \setminus \union[a \in T_\alpha][B_\alpha^a]$. As there are only countably
  many approximations, this set is \Borel. This completes the inductive
  construction.

  As there are only countably many approximations and $\sequence
  {S_\alpha}[\alpha < \omega_1]$ is decreasing, there exists $\alpha_0 <
  \omega_1$ such that $S_{\alpha_0 + 1} = S_{\alpha_0}$.

  \begin{claim} \label{ClaimDerivation}
    Every element of $S_{\alpha_0}$ has an immediate successor in
    $S_{\alpha_0}$.
  \end{claim}

  \begin{claimproof}
    Let $a$ be an approximation having no successor in $S_{\alpha_0}$.
    Then $a$ is $\alpha_0$-terminal, so $B_{\alpha_0}^a$ is defined and
    $\compatible{a}{X_{\alpha_0} \setminus B_{\alpha_0}^a} = \emptyset$.
    Since $X_{{\alpha_0} + 1} \subseteq X_{\alpha_0} \setminus
    B_{\alpha_0}^a$, it follows that $\compatible{a}{X_{{\alpha_0} + 1}} =
    \emptyset$, hence $a \notin S_{{\alpha_0} + 1}$. But $S_{\alpha_0 + 1}
    = S_{\alpha_0}$, so $\alpha \notin S_{\alpha_0}$.
  \end{claimproof}

  Let $a^0$ denote the unique approximation for which $n^{a^0} \defeq 0$
  and $\calU^{a^0} \defeq \continuousfunctions{\CantorCantorspace[0]}{X}$.
  There are two cases.

  \case{1}{$a^0 \notin S_{\alpha_0}$.} Then no configuration is
  simultaneously compatible with $a^0$ and generically compatible with
  $X_{\alpha_0}$. Observe that a configuration $\gamma$ such that
  $n^\gamma = 0$ essentially consists only of a continuous mapping
  $f^\gamma \from \CantorCantorspace[0] \to X$, and can therefore be
  identified with a point $x^\gamma \in X$. Every such configuration is
  compatible with $a^0$, and such a configuration is generically compatible
  with $X_{\alpha_0}$ if and only if $x^\gamma \in X_{\alpha_0}$. It follows
  that $X_{\alpha_0} = \emptyset$, so $X = \union[\alpha < \alpha_0]
  [X_\alpha \setminus X_{\alpha + 1}] = \union[\alpha < \alpha_0][{\union[a
  \in T_\alpha][B_\alpha^a]}]$. For all $n \in \N$, the set $B_n \defeq \union
  [\alpha < \alpha_0][{\union[a \in T_\alpha, n^a = n][B_\alpha^a]}]$ is a
  countable union of $E_n$-invariant \Borel sets on which $E$ has
  countable index over $E_n$, and is therefore itself an $E_n$-invariant
  \Borel set on which $E$ has countable index over $E_n$. As $X = \union
  [n \in \N][B_n]$, condition (1) follows.

  \case{2}{$a^0 \in S_{\alpha_0}$}. Then Claim \ref{ClaimDerivation} yields
  a sequence $\sequence{a^n}[n \in \N]$ of elements of $S_{\alpha_0}$
  such that $a^{n+1}$ is an immediate successor of $a^n$ for all $n \in \N$.
  It follows that $n^{a^n} = n$ for every $n \in \N$. As $a^n \in
  S_{\alpha_0}$, there is a configuration $\gamma^n$ which is compatible
  with $a^n$. For all $n \in \N$, define a continuous function $f^n \from
  \CantorCantorspace \to X$ by $f^n \defeq \sectionpower{(f^{\gamma^n})}
  {\N}$. Similarly, for all $n \in \N$ and $k < n$, define a continuous function
  $g_k^n \from \Fzero{k+1} \setminus \Fzero{k} \to Y_k$ by $g_k^n \defeq
  \sectionpower{(g_k^{\gamma^n})}{\N}$. If $m \le n$, then $f^{\gamma^m}
  \in \calU^{a^m}$ and $f^{\gamma^n} \in \calU^{a^n} \subseteq
  \sectionpower{(\calU^{a^m})}{n}$, so $f^m, f^n \in \sectionpower
  {(\calU^{a^m})}{\N}$. As our choice of metrics ensures that $\diam
  (\sectionpower{(\calU^{a^m})}{\N}) \le \diam(\calU^{a^m}) \le 1 / m$, it
  follows that $d(f^m, f^n) \le 1 / m$, so $\sequence{f^n}[n \in \N]$ is a
  \Cauchy sequence and admits a limit $f \in \continuousfunctions
  {\CantorCantorspace}{X}$. Similarly, for all $k \in \N$, the sequence
  $\sequence{g_k^n}[n > k]$ is \Cauchy, and therefore admits a limit $g_k
  \in \continuousfunctions{\Fzero{k+1} \setminus \Fzero{k}}{Y_k}$.

  We now show that $f$ satisfies condition (2). It is sufficient to show that if
  $k \in \N$ and $\pair{x}{x'} \in \Eone \setminus \Fzero{k}$, then $\pair{f(x)}
  {f(x')} \in E \setminus R_k$. As $\sequence{R_l}[l \in \N]$ is increasing
  and $\Eone \setminus \Fzero{k} = \union[l \ge k][\Fzero{l+1} \setminus
  \Fzero{l}]$, we can assume that $\pair{x}{x'} \in \Fzero{k+1} \setminus
  \Fzero{k}$. Observe that if $n > k$, then
  \begin{align*}
    \pair{f^n(x)}{f^n(x')}
      & = \pair{f^{\gamma^n}(x \restriction n)}{f^{\gamma^n}(x' \restriction n)}
        \\
      & = \pi_k(g_k^{\gamma^n}(x \restriction n, x' \restriction n)) \\
      & = \pi_k(g_k^n\pair{x}{x'}).
  \end{align*}
  By Proposition \ref{ContEval}, we can take the limit on both sides of the
  equality to obtain $\pair{f(x)}{f(x')} = \pi_k(g_k(x, x'))$, so $\pair{f(x)}{f(x')}
  \in E \setminus R_k$.
\end{theoremproof}

Before proving our second dichotomy, we note the following consequence
of \cite[Propositions 2.13 and 2.14]{TreeableEqRel}:

\begin{proposition} \label{kMeager}
  Let $C \subseteq \CantorCantorspace$ be a comeager set and
  $\sequence{R_n}[n \in \N]$ be a sequence of \Fsigma binary relations
  on $\CantorCantorspace$ such that $R_n \intersection \Eone
  \setminus \Fzero{n} = \emptyset$ for all $n \in \N$. Then there is a
  continuous homomorphism $\phi \from \CantorCantorspace \to C$ from
  $\sequence{\Fzero{n}, \setcomplement{\Fzero{n}}}[n \in \N]$ to
  $\sequence{\Fzero{n}, \setcomplement{R_n}}[n \in \N]$.
\end{proposition}

Our second dichotomy is the following:

\begin{theorem} \label{UnionReduction}
  Let $X$ be a \Polish space, $E$ be an analytic equivalence relation on
  $X$, and $\sequence{E_n}[n \in \N]$ be a sequence of potentially
  \Fsigma subequivalence relations of $E$. Then exactly one of the
  following holds:
  \begin{enumerate}
    \item There is a cover $\sequence{B_n}[n \in \N]$ of $X$ such that
      $B_n$ is an $E_n$-invariant \Borel set on which $E$ has countable
      index over $E_n$ for all $n \in \N$.
    \item There is a continuous homomorphism $\phi \from
      \CantorCantorspace \to X$ from $\brokenpair{\Eone}{\sequence
      {\setcomplement{\Fzero{n}}}[n \in \N]}$ to $\pair{E}{\sequence
      {\setcomplement{\union[i \le n][E_i]}}[n \in \N]}$.
  \end{enumerate}
  Moreover, if $E = \union[n \in \N][E_n]$, then the mapping $\phi$ in the
  latter condition is an embedding of $\Eone$ into $E$.
\end{theorem}

\begin{theoremproof}
  The ``moreover'' part is immediate. Condition (2) in the statement of the
  theorem is stronger than condition (2) in the statement of Theorem \ref
  {UnionCohomomorphism}; hence, as in Theorem \ref
  {UnionCohomomorphism}, the two conditions are mutually exclusive. It
  remains to show that at least one of them holds.

  Refining the topology on $X$ if necessary, we can assume that the
  $E_n$'s are \Fsigma. We need to show that if condition (2) in the
  statement of Theorem \ref{UnionCohomomorphism} holds, then so too
  does condition (2) in the statement of Theorem \ref{UnionReduction}.
  Suppose that $\psi \from \CantorCantorspace \to X$ is a continuous
  homomorphism from $\sequence{\Eone \setminus \Fzero{n}}[n \in \N]$ to
  $\sequence{E \setminus \union[i \le n][E_i]}[n \in \N]$. For all $n \in \N$,
  define $R_n' \defeq \preimage{(\psi \times \psi)}{\union[i \le n][E_i]}$. Then
  $R_n'$ is \Fsigma and $R_n' \intersection (\Eone \setminus \Fzero{n}) =
  \emptyset$. So Proposition \ref{kMeager} yields a continuous
  homomorphism $\pi \from \CantorCantorspace \to \CantorCantorspace$
  from $\sequence{\Fzero{n}, \setcomplement{\Fzero{n}}}[n \in \N]$ to
  $\sequence{\Fzero{n}, \setcomplement{R_n'}}[n \in \N]$, in which cases
  the mapping $\phi \defeq \psi \composition \pi$ satisfies condition (2).
\end{theoremproof}

Following the usual abuse of language, we say that an equivalence relation
is \definedterm{finite} if each of its equivalence classes are finite. A \Borel
equivalence relation $E$ on a \Polish space is \definedterm{hyperfinite} if
there is an increasing sequence $\sequence{E_n}[n \in \N]$ of finite \Borel
subequivalence relations of $E$ whose union is $E$. As the
\Lusin--\Novikov uniformization theorem easily implies that every finite
\Borel equivalence relation on a \Polish space is smooth, it follows that
every hyperfinite \Borel equivalence relation on a \Polish space is
hypersmooth.

We will need the following elementary result (see \cite[Proposition 5.1]
{TreeableEqRel}):

\begin{proposition} \label{EssHyperfinite}
  Suppose that $X$ is a \Polish space, $\sequence{E_n}[n \in \N]$ is an
  increasing sequence of smooth \Borel equivalence relations on $X$, and
  there is a cover $\sequence{B_n}[n \in \N]$ of $X$ such that $B_n$ is a
  \Borel set on which $\union[m \in \N][E_m]$ has countable index over
  $E_n$ for all $n \in \N$. Then $\union[m \in \N][E_m]$ is \Borel reducible to
  a hyperfinite \Borel equivalence relation on a \Polish space.
\end{proposition}

As a first application of our second dichotomy, we give a proof of
the \Kechris--\Louveau dichotomy:

\begin{theorem}[\Kechris--\Louveau] \label{KechrisLouveau}
  Suppose that $E$ is a hypersmooth \Borel equivalence relation
  on a \Polish space. Then exactly one of the following holds:
  \begin{enumerate}
    \item There is a \Borel reduction of $E$ to a hyperfinite \Borel
      equivalence relation on a \Polish space.
    \item There is a continuous embedding of $\Eone$ into $E$.
  \end{enumerate}\end{theorem}

\begin{theoremproof}
  The exclusivity of the two conditions comes from the fact, mentioned in
  the introduction, that $\Eone$ is not \Borel reducible to a countable \Borel
  equivalence relation on a \Polish space; this can, for instance, be
  obtained as a consequence of \cite[Theorem 4.1]{KechrisLouveau}, \cite
  [\S1.II.i]{Kechris:LocallyCompact}, and \Feldman--\Moore's theorem
  (alternatively, a more elementary proof can be obtained from \cite
  [Proposition 2.4]{TreeableEqRel} and \cite[Proposition 2.5]
  {TreeableEqRel}). To see that at least one of the two conditions holds,
  write $E$ as the increasing union of a sequence of smooth \Borel
  subequivalence relations $\sequence{E_n}[n \in \N]$. Refining the
  topology on $X$ if necessary, we can assume that the $E_n$'s are closed.
  Hence we can apply Theorem \ref{UnionReduction} and use Proposition
  \ref{EssHyperfinite} to complete the proof.
\end{theoremproof}

\section{Primary results} \label{finalSection}

The following fact is the main result of this paper:

\begin{theorem} \label{finalSection:main}
  Suppose that $\calF$ is a class of strongly-idealistic potential\-ly-\Fsigma
  equivalence relations on \Polish spaces that is closed under countable
  disjoint unions and countable-index \Borel superequivalence relations.
  If $E$ is an equivalence relation on a \Polish space that is a countable
  union of subequivalence relations that are \Borel reducible to relations in
  $\calF$, then at least one of the following holds:
  \begin{enumerate}
    \item There is a \Borel reduction of $E$ to a relation in $\calF$.
    \item There is a continuous embedding of $\Eone$ into $E$.
  \end{enumerate}
  Moreover, if every relation in $\calF$ is also ccc idealistic, then exactly
  one of these conditions holds.
\end{theorem}

\begin{theoremproof}
  To see that the conditions are mutually exclusive when every relation in
  $\calF$ is ccc idealistic, observe that if both hold, then there is a \Borel
  reduction of $\Eone$ to a ccc idealistic \Borel equivalence relation on a
  \Polish space, contradicting \cite[Theorem 4.1]{KechrisLouveau}.

  To see that at least one of the conditions holds, note that, by Theorem
  \ref{UnionReduction}, we can assume that there is a cover $\sequence
  {B_n}[n \in \N]$ of $X$ by \Borel sets on which $E$ has countable index
  over subequivalence relations that are \Borel reducible to relations in
  $\calF$. For each $n \in \N$, Proposition \ref{SwitchSidesCountableIndex}
  yields a \Borel reduction $\phi_n$ of $E \restriction B_n$ to some $F_n
  \in \calF$, and Proposition \ref{InvariantReduction} allows us to assume
  that $B_n$ is $E$-invariant. Then the $E$-invariant \Borel sets $B_n'
  \defeq B_n \setminus \union[m < n][B_m]$ partition $X$, so the functions
  $\phi_n \restriction B_n'$ can be combined to obtain a \Borel reduction of
  $E$ to $\disjointunion[n \in \N][F_n]$.
\end{theoremproof}

Theorem \ref{introduction:main} follows from Theorem \ref
{finalSection:main} and the previously mentioned fact that every countable
\Borel equivalence relation on a \Polish space is strongly ccc idealistic and
potentially \Fsigma.

Theorem \ref{introduction:action} follows from Theorem \ref
{introduction:main}, \cite[Theorem 4.1]{KechrisLouveau}, the
\Feldman--\Moore theorem, and the fact that every countable \Borel
equivalence relation on a \Polish space is ccc idealistic.

Theorem \ref{introduction:idealistic} follows from Proposition
\ref{idealistic:countableindex}, Theorem \ref{finalSection:main}, and
\cite[Theorem 4.1]{KechrisLouveau}.

Theorem \ref{introduction:unions} follows from Theorem \ref
{introduction:unions:two} and the fact that every countable \Borel
equivalence relation on a \Polish space is strongly ccc idealistic and
potentially \Fsigma.

Theorem \ref{introduction:unions:two} follows from Theorem \ref
{finalSection:main} and \cite[Theorem 4.1]{KechrisLouveau}.

Theorem \ref{introduction:intersecting} follows from Theorem \ref
{introduction:unions:two}, Propositions \ref
{CountableIndexToIntersectionReduction}, \ref
{IntersectingCountableUnion}, and \ref{idealistic:countableindex}, and the
observation that if $F$ is a \Borel equivalence relation on a \Polish space,
then the class $\calF$ of equivalence relations on \Polish spaces that are
\Borel isomorphic to countable-index \Borel superequivalence relations of
$F \times \diagonal{\N}$ is closed under countable disjoint unions and
countable-index \Borel superequivalence relations, and if $F$ is strongly
idealistic and potentially \Fsigma, then so too is $F \times \diagonal
{\N}$.

\begin{acknowledgements}
  We would like to thank Zolt\'an Vidny\'anszky and Jind\v{r}ich Zapletal for
  several enlightening discussions.
\end{acknowledgements}

\bibliographystyle{amsalpha}
\bibliography{main}

\end{document}